\documentclass[12pt]{article}
\usepackage{fancyhdr}
\usepackage{fullpage}
\usepackage{amsmath}
\usepackage{amssymb}
\usepackage{theorem}
\newtheorem{thm}{Theorem}[section]
\newtheorem{prop}{Proposition}[section]
\newtheorem{corr}{Corollary}[section]
\newtheorem{lemma}{Lemma}[section]
{\theorembodyfont{\rm} \newtheorem{exam}{Example}[section]
\newtheorem{rem}{Remark}[section]}
\newtheorem{defe}{Definition}[section]
\newtheorem{conj}{Conjecture}[section]
\numberwithin{equation}{section}
\def\sab{\text{sign}(\alpha, \beta)}
\def\kab{K_{\alpha,\beta}}
\def\pab{(\alpha, \beta)}
\def\qh{(q - q^{-1})}
\def\eab{e_{\beta} \otimes e_{-\alpha}}

\def\qs{$\quad \square$}
\def\R{\mathbb R}

\def\Res{\text{Res}}

\def\ell{{\text{ell}}}

\def\h1{\hat{\bold 1}}

\def\Ua{U_q(\tilde\g)}
\def\U2{{\Ua}_2}
\def\g{\mathfrak g}

\def\N{\mathbb N}
\def\Z{\mathbb Z}
\def\C{\mathbb C}

\def\<{\langle}
\def\>{\rangle}
\def\o{\otimes}

\def\h{{\mathfrak h}}
\def\g{{\mathfrak g}}

\begin{document}
\thispagestyle{plain}
\begin{center}
{\Large \bf On the GGS Conjecture}

\vskip 12 pt

{\large \bf Travis Schedler}
\end{center}
\begin{abstract}
In the 1980's, Belavin and Drinfeld classified solutions $r$ of the
classical Yang-Baxter equation (CYBE) for simple Lie algebras
$\mathfrak g$ satisfying $0 \neq r + r_{21} \in (S^2
\mathfrak{g})^{\mathfrak{g}}$ \cite{BD}.  They proved that all such
solutions fall into finitely many continuous families and introduced
combinatorial objects to label these families, Belavin-Drinfeld
triples. In 1993, Gerstenhaber, Giaquinto, and Schack attempted to
quantize such solutions for Lie algebras $\mathfrak{sl}(n).$ As a
result, they formulated a conjecture stating that certain explicitly
given elements $R \in Mat_n(\mathbb C) \otimes Mat_n(\mathbb C)$
satisfy the quantum Yang-Baxter equation (QYBE) and the Hecke relation
\cite{GGS}.  Specifically, the conjecture assigns a family of such
elements $R$ to any Belavin-Drinfeld triple of type $A_{n-1}$.
Following a suggestion from Gerstenhaber and Giaquinto, we propose an
alternate form for $R$, given by $R_J = q^{r^0} J^{-1} R_s J_{21}
q^{r^0}$, for a suitable twist $J$ and a diagonal matrix $r^0$, where
$R_s$ is the standard Drinfeld-Jimbo solution of the QYBE.  We
formulate the ``twist conjecture'', which states that $R_J =
R_{\text{GGS}}$ and that $R_J$ satisfies the QYBE.  Since $R_J$ by
construction satisfies the Hecke relation, this conjecture implies the
GGS conjecture. We check the twist conjecture by computer for $n \leq
12$ and show that it is true modulo $\hbar^3$.  We provide
combinatorial formulas for coefficients in the matrices $R_J,
R_{\text{GGS}}$ and prove both conjectures in the {\it orthogonal
generalized disjoint case}---where $\Gamma_1 = \bigcup_i \Gamma_1^{i}$
with $\Gamma_1^{i} \perp \Gamma_1^{j}, i \neq j$, $\tau \Gamma_1^{i}
\cap \Gamma_1 \subset \Gamma_1^{i+1}$, and $\tau^j \Gamma_1^{i} \perp
\Gamma_1^{i}, \forall i,j \geq 1$. We also prove the twist conjecture
in the {\it disjoint case}, $\Gamma_1 \cap \Gamma_2 = \emptyset$.
Finally, we prove the twist conjecture for the Cremmer-Gervais triple
and discuss cases in which it is known that $R_J = R_{\text{GGS}}$.
\end{abstract}

\section{Main Results}
We begin this section by introducing Belavin-Drinfeld triples.  We
present the GGS conjecture, which is motivated by calculating possible
quantizations of $r$ modulo $\hbar^3$.  Next, we proceed to formulate
the twist conjecture and give the remarkably similar combinatorial
descriptions of the twist and the GGS $R$-matrix.  Finally, we summarize our
main results, namely the computer verification of the twist
conjecture, its proof modulo $\hbar^3$, and a complete proof of the
twist conjecture in the disjoint, orthogonal generalized disjoint, and
Cremmer-Gervais cases.
%  We provide a combinatorial
%description and of the coefficients of $R_{\text{GGS}} and J$.  Then,
%we summarize our main results, namely the computer verification of the
%conjectures, their equivalence modulo $\hbar^3$, and their proof in
%the disjoint, orthogonal generalized disjoint, and Cremmer-Gervais
%cases.
\subsection{Belavin-Drinfeld triples}

Let $(e_i), 1 \leq i \leq n,$ be a basis for $\mathbb C^n$.  Set
$\Gamma = \{e_i - e_{i+1}: 1 \leq i \leq n-1\}$.  We will use the
notation $\alpha_i \equiv e_i - e_{i+1}$.  Let $( , )$ denote the
inner product on $\mathbb C^n$ having $(e_i)$ as an orthonormal basis.

\begin{defe} \cite{BD}
A Belavin-Drinfeld triple of type $A_{n-1}$ is a triple 
$(\tau, \Gamma_1, \Gamma_2)$ where \\ $\Gamma_1, \Gamma_2 \subset \Gamma$
and $\tau: \Gamma_1 \rightarrow \Gamma_2$ is a bijection, satisfying
two conditions:

(a) $\forall \alpha, \beta
\in \Gamma_1$, $(\tau \alpha,\tau \beta) = (\alpha, \beta)$.

(b) $\tau$ is nilpotent: $\forall \alpha \in \Gamma_1, \exists k
\in \mathbb N$ such that $\tau^k \alpha \notin \Gamma_1$.
\end{defe}

Let $\mathfrak g = \mathfrak{gl}(n)$ be the Lie algebra of $n \times
n$ matrices. (Although $\mathfrak{gl}(n)$ is not simple, solutions
correspond to those in $\mathfrak{sl}(n)$, and it will simplify
computations.  For the same reason, we state the GGS and twist
conjectures in $\mathfrak{gl}(n)$.)  Set $\mathfrak h \subset
\mathfrak g$ to be the subset of diagonal matrices.  Elements of
$\mathbb C^n$ define linear functions on $\mathfrak h$ by $\bigl
( \sum_i \lambda_i e_i \bigr) \bigl( \sum_i a_i\: e_{ii} \bigr)=
\sum_i \lambda_i a_i$.  Let $P = \sum_{1 \leq i,j \leq n} e_{ij}
\otimes e_{ji}$ be the Casimir element for $\mathfrak g$ as well as
the permutation matrix, and let $P^0=\sum_i e_{ii} \o e_{ii}$ be the
projection of $P$ to $\mathfrak h \otimes \mathfrak h$.
 
For any Belavin-Drinfeld triple, consider the following equation for
$r^0 \in \mathfrak h \wedge \mathfrak h$:

\begin{gather} 
%\label{tr01}
%r^0_{12} = - r^0_{21}. \\ 
\label{tr02} \forall \alpha \in \Gamma_1,
\bigl[(\alpha - \tau \alpha) \otimes 1 \bigr]r^0 = \frac{1}{2} 
\bigl[(\alpha + \tau \alpha) \otimes 1\bigr] P^0.
\end{gather}

Belavin and Drinfeld showed that solutions $r \in \mathfrak{g} \o
\mathfrak g$ of the CYBE satisfying $r + r^{21} = P$, up to
isomorphism, are given by a discrete datum (the Belavin-Drinfeld
triple) and a continuous datum (a solution $r^0 \in \h \wedge \h$ of
\eqref{tr02}).  We now describe this classification.  For $\alpha =
e_i - e_j$, set $e_\alpha \equiv e_{ij}$, and say $\alpha > 0$ if $i <
j$, and otherwise $\alpha < 0$.  Define $|\alpha| = |j - i|$.  For any
$Y \subset \Gamma$, set $\tilde Y = \{v \in \text{Span}(Y) \mid v =
e_i - e_j, v > 0\}$; in particular we will use $\tilde \Gamma_1,
\tilde \Gamma_2$.  We extend $\tau$ additively to a map $\tilde
\Gamma_1 \rightarrow \tilde \Gamma_2$, i.e.  $\tau(\alpha+\beta)=\tau
\alpha +\tau \beta$. Whenever $\tau^k \alpha = \beta$ for $k \geq 1$,
we say $\alpha \prec \beta$.  Clearly $\prec$ is a partial ordering on
$\tilde \Gamma$.  Finally, for any $\beta = \tau^k \alpha$, $\alpha=
e_j - e_i, \beta = e_k - e_{k+i-j}, j < i-1,$ we say $\tau^k$ reverses
orientation on $\alpha$ if $\tau \alpha_j = \alpha_{k+i-j-1}$ and
$\tau^k$ preserves orientation on $\alpha$ if $\tau \alpha_j =
\alpha_k$.  In the reversing case, write $\alpha \prec^{\leftarrow}
\beta$ and $\text{sign}(\alpha,\beta) = (-1)^{1-|\alpha|}$; in the
preserving case write $\alpha \prec^\rightarrow \beta,
\text{sign}(\alpha,\beta) = 1$. We set $\text{sign}(\alpha,\beta) = 1$
when $j+1=i$.  Set $x \wedge y \equiv x \o y - y \o x$ for $x,y \in
Mat_n(\C)$ and $z = \sum_{i,j,k,l} z_{ik}^{jl} e_{ij} \o e_{kl}$ for
$z \in Mat_n(\C) \o Mat_n(\C)$. Then we define

\begin{equation}
a = \sum_{\alpha \prec \beta} \text{sign}(\alpha,\beta)\: e_{-\alpha}
\wedge e_{\beta}, \quad r_s = \frac{1}{2} \sum_i e_{ii} \o e_{ii} +
\sum_{\alpha > 0} e_{-\alpha} \o e_{\alpha}, \quad r = r^0 + a +
r_s, \label{r}
\end{equation}
where $r_s \in \mathfrak{g} \o \mathfrak{g}$ is the standard solution
of the CYBE satisfying $r_s + r_s^{21} = P$, and $r$ is the solution
corresponding to the data $((\Gamma_1,\Gamma_2,\tau), r^0)$.  It
follows from \cite{BD} that any solution $\tilde r \in \mathfrak{g},
\tilde r+\tilde r_{21} = P$ is equivalent to such a solution $r$ under
an automorphism of $\mathfrak{g}$.

\subsection{The GGS conjecture}
The GGS conjecture proposes a hypothetical quantization of the matrix
$r$ given in \eqref{r}, given by a matrix $R \in Mat_n(\C) \o
Mat_n(\C)$ conjectured to satisfy the quantum Yang-Baxter equation
(QYBE), $R_{12} R_{13} R_{23} = R_{23} R_{13} R_{12}$, and the Hecke
relation, $(PR - q)(PR+q^{-1}) = 0$.  This may be formulated and justified
as follows (which is more or less the original motivation).

If we write $R \equiv 1 + 2 \hbar r + 4 \hbar^2 s
\pmod{\hbar^3}$, where $q \equiv e^\hbar$, then we can consider the constraints
imposed by the QYBE and the Hecke relation modulo $\hbar^3$.  One may
easily check that QYBE becomes the CYBE for $r$, while the Hecke
relation becomes the condition $s + s_{21} = r^2$.  Thus, there is a
unique choice of $s$ that is symmetric, namely $\frac{1}{2} r^2 =
\frac{1}{2} ((r^0)^2 + a r^0 + r^0 a + \epsilon)$ where
\begin{equation} \label{eps}
\epsilon = ar_s + r_s a + a^2.
\end{equation}
%Set $c = \sum_{\alpha > 0} e_{-\alpha} \wedge e_\alpha$; then
%$a r_s + r_s a = ac + ca$ since $a$ is skew-symmetric.  

\begin{prop} \label{ggs1}
There exist unique polynomials $P_{i,j,k,l}$ of the form
$x q^y (q-q^{-1})^z, x,y \in \C, z \in \{0,1\}$ such that 
$\sum_{i,j,k,l} P_{i,j,k,l} e_{ij} \o e_{kl} \equiv 1 + 2\hbar r
+2 \hbar^2 r^2 \pmod{\hbar^3}$.
\end{prop}

{\it Proof.}  The proof is easy.   \qs

\begin{defe} Define $R_{\text{GGS}} = \sum_{i,j,k,l} P_{i,j,k,l} e_{ij} \o e_{kl}$, with the $P_{i,j,k,l}$ uniquely determined by Proposition \ref{ggs1}.  The
matrix $R_{\text{GGS}}$ is called the GGS R-matrix.
\end{defe}

Define the following matrices:

\begin{equation}
\tilde a = \sum_{i,j,k,l} a_{ik}^{jl} q^{a_{ik}^{jl} \epsilon_{ik}^{jl}}
e_{ij} \o e_{kl},
\quad \bar R_{\text{GGS}} = R_s + (q-q^{-1}) \tilde a,
\end{equation}
where 
$R_s = q \sum_{i} e_{ii} \otimes e_{ii} +
\sum_{i \neq j} e_{ii} \otimes e_{jj} + (q - q^{-1}) \sum_{i>j} e_{ij}
\otimes e_{ji}$ is the standard Drinfeld-Jimbo solution to the QYBE.

\begin{prop}  The matrix $R_{\text{GGS}} = q^{r^0} \bar R_{\text{GGS}} 
q^{r^0}$.
\end{prop}

{\it Proof.} This is clear. \qs

\begin{rem} We see that $R_{\text{GGS}} \equiv q^{2r} \pmod{\hbar^3}$,
although $R_{\text{GGS}} \neq q^{2r}$ in general.
\end{rem}

\begin{conj}{\bf ``the GGS conjecture'' \cite{GGS}} \label{ggs}  

I. The matrix $R_{\text{GGS}}$ satisfies the QYBE. 

II. The matrix $R_{\text{GGS}}$ satisfies the Hecke relation.
\end{conj}

We will sometimes refer separately to the two parts as Conjectures \ref{ggs}.I
and \ref{ggs}.II.  

\begin{rem} \label{or0} It is sufficient to check the QYBE for one $r^0$
since the space of homogeneous solutions to \eqref{tr02} is exactly the
space $\Lambda^2(\mathfrak{l})$ where $\mathfrak{l} \subset \mathfrak{h}$
is the space of symmetries of the Belavin-Drinfeld triple, i.e. $(x,\alpha)
= (x, \tau \alpha)$ for any $x \in \mathfrak{l}, \alpha \in \Gamma_1$.
It is easy to see that $x \in \mathfrak{l}$ implies $[1 \o x + x \o 1,
R_{\text{GGS}}] = 0$, and it follows for any $y \in \Lambda^2(\mathfrak{l})$
that $q^y R_{\text{GGS}} q^y$ satisfies the QYBE iff $R_{\text{GGS}}$ does. 
The same holds for $R_J$ as defined in the following section.
\end{rem}

Now, we describe our new results on the GGS conjecture.

\begin{thm} \label{ggscomp}
(i) The GGS conjecture is true for $n \leq 12$.
(ii) The GGS conjecture is true modulo $\hbar^3$.
\end{thm}

{\it Proof.} (i) This has been verified by the author through computer programs
detailed in \cite{S}.  The programs check the QYBE and Hecke relation
directly using one choice of $r^0$.  One may check that this is
sufficient to prove GGS for any $r^0$.

(ii) This is obvious from construction.
\qs

We see that the strangest matrix in the definition of $R_{\text{GGS}}$
is $\epsilon$.  Here we give a simple combinatorial formula for this
unusual matrix.  For $i<j, k < l$, we say that $e_i - e_j \lessdot e_k
- e_l$ iff $j = k$, and similarly define $\gtrdot$.  Let
$[\text{statement}]=1$ if ``statement'' is true and $0$ if
``statement'' is false.

\begin{prop}  \label{fep} We may rewrite $\epsilon$ as follows:
\begin{multline} \label{fe}
\epsilon = \sum_{\alpha \prec \beta} \text{sign}(\alpha,\beta)\bigl[
-\frac{1}{2} [\alpha \lessdot \beta] - \frac{1}{2}[\beta \lessdot
\alpha] - [\exists \gamma, \alpha \prec \gamma \prec \beta, \alpha
\lessdot \gamma] \\ - [\exists \gamma, \alpha \prec \gamma \prec
\beta, \alpha \gtrdot \gamma] + [\alpha \prec^{\leftarrow} \beta] (1 -
|\alpha|) \bigr] (\eab + e_{-\alpha} \o e_{\beta})
\end{multline}
\end{prop}

{\it Proof.} Given in Section 2. \qs

\begin{exam} \label{gcge} For a given $n$, there
are exactly $\phi(n)$ triples ($\phi$ is the Euler $\phi$-function) in
which $|\Gamma_1| + 1 = |\Gamma|$ \cite{GG}.  These are are called
{\it generalized Cremmer-Gervais} triples.  These are indexed by $m
\in \Z^+$, where $\text{gcd}(n,m) = 1$, and given by $\Gamma_1 =
\Gamma \setminus \{\alpha_{n-m}\}$, $\Gamma_2 = \Gamma \setminus
\{\alpha_m\}$, and $\tau(\alpha_i) = \alpha_{\text{Res}(i+m)}$, where
$\Res$ gives the residue modulo $n$ in $\{1,\ldots,n\}$.  For these
triples, there is a unique $r^0$ with first component having trace 0,
which is given by $(r^0)^{ii}_{ii} = 0, \forall i$, and
$(r^0)_{ij}^{ij} = \frac{1}{2} - \frac{1}{n}\text{Res}(\frac{j-1}{m})$
(this is easy to verify directly and is also given in \cite{GG}).
With this $r^0$, $R_{\text{GGS}}$ has a very nice combinatorial
formula, which was conjectured by Giaquinto and checked in some cases.
We now state and prove this formula. As in \cite{GH}, define
$e_{-\alpha} \wedge_c e_\beta = q^{-c} e_{-\alpha} \otimes e_{\beta} -
q^{c} e_{\beta} \otimes e_{-\alpha}$. Let $O(\alpha,\beta) = l$
when $\tau^l \alpha = \beta$.
\begin{prop} $R_{\text{GGS}}$ is given as follows:
\begin{equation} \label{gcgr}
R_{\text{GGS}} = q^{r^0} R_s q^{r^0} + \sum_{\alpha \prec \beta} (q -
q^{-1}) e_{-\alpha} \wedge_{\frac{-2O(\alpha,\beta)}{n}} e_\beta.
\end{equation} 
\end{prop}
{\it Proof}. See Appendix B. \qs
\end{exam}

\begin{rem}  Our formulation is from \cite{GH}, correcting
misprints.  The original formulation in \cite{GGS} is somewhat
different. We will write $x_{q^{-1}}$ to denote the matrix $x$ with
$q^{-1}$ substituted for $q$. Define $(x \o y)^T = x^T \o y^T$ where
$x^T$ is the transpose of $x$, for $x, y \in Mat_n(\C)$. Then, the
original form of $R_{\text{GGS}}$ can be written as follows:
$$R = q^{-r^0} \bigl( R_s + (q^{-1} - q) \tilde a_{q^{-1}}^T
\bigr) q^{-r^0}.$$ 
Denoting $R$ as this matrix and
$R_{\text{GGS}}$ as given before, we have $R_{\text{GGS}} -
R^T_{q^{-1}} = q^{r^0} (q - q^{-1}) P q^{r^0}$ \\ $= (q -
q^{-1}) P$.  Thus, $R_{\text{GGS}}$ satisfies the Hecke relation iff
$R$ satisfies the Hecke relation.  In this case, we have $P
R^T_{q^{-1}} = (P R_{\text{GGS}})^{-1}$, so $R^T_{q^{-1}} = 
(R_{\text{GGS}}^{-1})_{21}$, and thus $R$ satisfies the QYBE iff
$R_{\text{GGS}}$ does.  Thus, the two formulations are equivalent.
\end{rem}

\subsection{The twist conjecture} \label{j1}

%There is a problem with the GGS conjecture: it is not obvious why such
%a formula should give a solution of the Yang-Baxter equation, and the
%matrix $\epsilon$ is particularly mystifying.  To address this issue,
%two steps have been taken.  First, since the conjecture is so simple
%to state, it lends well to computer verification, which has been
%carried out to Lie algebras $\mathfrak{sl}(n)$ for $n \leq 12$ by the
%author, as detailed in \cite{S}.  Here we take a second step: we
%attempt to explain why the conjecture makes sense by assuming that
%$R_{\text{GGS}}$ may be rewritten in a form $R_J = q^{r^0}
%J^{-1} R_s J_{21} q^{r^0}$ where $R_J = \prod_{\alpha \prec
%\beta} (1 + (q-q^{-1}) q^{K_{\alpha,\beta}} \eab)$ for a suitable
%ordering on pairs $\alpha \prec \beta$.  We find that under this
%assumption, $\epsilon$ is uniquely determined if assumed to be
%symmetric.  Furthermore, we shed light on $\epsilon$ and
%$K_{\alpha,\beta}$ by giving straightforward combinatorial formulas.

In \cite{EK}, it is proved that any quasitriangular structure as
defined in Section 1.1 has a quantization which is a twist of the
standard quasitrangular Hopf algebra $U_q(\mathfrak{gl}(n))$.  In
\cite{H} (see also \cite{ER}), such a twist is constructed for the
{\it disjoint} case, $\Gamma_1 \cap \Gamma_2 = \emptyset$ (another
twist is given in Appendix \ref{gda}).  Thus, in the $n$-dimensional
representation, there should exist $J \in Mat_n(\C) \o Mat_n(\C)$ so
that $R_J = q^{r^0} J^{-1} R_s J_{21} q^{r^0}$ satisfies the QYBE and
Hecke relation.  Further, it is especially nice to look for {\it
triangular} twists: twists where $J = 1 + N$ and $N = \sum_{\alpha,
\beta \in \tilde \Gamma} N_{\alpha,\beta} \eab$.  In this section we
give an explicit $J$ of this form designed so that $R_J \equiv
R_{\text{GGS}} \pmod{\hbar^3}$ and conjecture that $R_J$ satisfies the
QYBE and $R_J = R_{\text{GGS}}$.  To find this interesting twist, the
author found the Gauss decomposition $\bar R_{\text{GGS}} = J^{-1} R_s
J_{21}$ (which is necessarily unique, if it exists) for all triples $n
\leq 12$.

First, we will define some useful notation.  Given a matrix
$$x = \sum_{\alpha, \beta > 0} N^+(\alpha, \beta) e_{\beta}
\otimes e_{-\alpha} + \sum_{\alpha, \beta > 0} N^-(\alpha, \beta)
e_{-\alpha} \otimes e_{\beta} + D,$$
where $D = \sum_i D_i \otimes D'_i$ with $D_i, D'_i$ diagonal, denote
\begin{gather*}
x_+ \equiv \sum_{\alpha, \beta > 0} N^+(\alpha, \beta) e_{\beta}
\otimes e_{-\alpha} + \frac{1}{2} D, \quad x_- \equiv \sum_{\alpha,
\beta > 0} N^-(\alpha, \beta) e_{-\alpha} \otimes e_{\beta} +
\frac{1}{2} D, \\ x_{\alpha,\beta} = N^+(\alpha, \beta), \quad
x_{-\beta,-\alpha} = N^-(\alpha,\beta).
\end{gather*}

Now, we proceed to define $J$.  Set $X = \{ (\alpha, \beta) \in \tilde
 \Gamma_1 \times \tilde \Gamma_2 \mid \alpha \prec \beta \}$ and $X^i
 = \{ (\alpha, \beta) \in X \mid \tau^i(\alpha) = \beta \}$ so that $X
 = \cup X^i$. Given any total ordering $<$ on a set $Y$, we will use
 $\displaystyle \sideset{}{^<}\prod_{x \in Y}$ to denote a product
 over all elements of $Y$, left to right, under the order $<$.

Define the following matrices, products taken left to right, 
with $K_{\alpha, \beta} \in \C$:
\begin{gather} 
\label{taid}
A^{i} = (q - q^{-1}) \sum_{{\beta = \tau^i(\alpha)}} \text{sign}(\alpha,
\beta) q^{K_{\alpha, \beta}} e_{\beta} \otimes e_{-\alpha}\\ \label{jd}
J^{i} = 1 + A^{i}, \quad J = \prod_{i = 1}^d
J^{i}, \quad \bar R_J = J^{-1} R_s
J_{21}, \quad R_J = q^{r^0} \bar R_J q^{r^0}.
\end{gather}

\begin{prop} There exists an ordering $<$ on $X$,
such that $J$ and $J^{-1}$ are given by the formulas
\begin{gather} \label{jprod}
J = \sideset{}{^<}\prod_{(\alpha, \beta) \in X} \bigl( 1 + \sab \qh q^{\kab} \eab \bigr), \\
\label{jiprod} J^{-1} = \sideset{}{^>}\prod_{(\alpha, \beta) \in X} \bigl( 1 - 
\sab \qh q^{\kab} \eab \bigr).
\end{gather}
\end{prop}

{\it Proof}.  Indeed, each $X^p$ may be ordered as follows: set $\beta
= e_i - e_j, \beta' = e_k - e_l$.  Then if $i > k$, $(\alpha, \beta) <
(\alpha', \beta')$.  If $i = k$ and $j > l$ then $(\alpha, \beta) <
(\alpha', \beta')$.  Then, it is clear that $\sideset{}{^<}\prod_{\pab
\in X^i} \bigl( 1 + \sab \qh q^{\kab} \eab \bigr) = J^i$ because, upon
expansion, all products $(\eab) (e_{\beta'} \otimes e_{-\alpha'})$
vanish.  Then, all that is needed is to extend $<$ to an ordering on
$X$ given by $X_i < X_j$ whenever $i < j$.  $\quad \square$

\begin{rem}
The product formula \eqref{jprod} for $J$ is especially natural in
light of the formula for the universal $R$-matrix given in \cite{KhT}.
In Appendix A, the importance of the ordering by powers of $\tau$ is
demonstrated in the construction of a twist ${\cal J} \in
U_q(\mathfrak{gl}(n)) \o U_q(\mathfrak{gl}(n))$ corresponding to $J$
for the simplest case where $\tau^2$ is defined (i.e. $\Gamma_1 \cap
\Gamma_2 \neq \emptyset$).  This is the orthogonal generalized
disjoint case.
\end{rem}
%\begin{rem}
%In Section \ref{ds}, we will see that, when
%$\Gamma_1 \cap \Gamma_2 = \emptyset$, the products \eqref{jprod},
%\eqref{jiprod} turn into sums of the same kind.  It is worth noting,
%however, that in general there is no such sum formulation of $J$.
%Indeed, already in the simple triple $n = 6$, $\Gamma_1 = \{\alpha_1,
%\alpha_2, \alpha_3\}$, $\tau \alpha_i = \alpha_{i+2}$, we obtain
%$J_{\alpha_1+\alpha_2,\alpha_4+\alpha_5} = (q - q^{-1})^2$, while 1)
%$\alpha_1+\alpha_2 \nprec \alpha_4+\alpha_5$, and 2) $(q - q^{-1})^2$
%is not of the form $\pm q^k (q - q^{-1}), k \in \C$.  For $n > 6$, one
%obtains terms not even of the form $m q^k(q-q^{-1})^l$, $m,k,l \in
%\C$.
%\end{rem}

%Set $q \equiv e^{\hbar}$.  Then we have the following result:

\begin{thm} \label{nkt}
%(i) For any choice of $\kab$, $\frac{d^2}{d\hbar^2} \bigl|_{\hbar=0}
%\bigl[\bar R_{J} + (\bar R_J)_{21}\bigr] =
%\frac{d^2}{d\hbar^2} \bigl|_{\hbar=0} \bigl[\bar R_{\text{GGS}}+ (\bar
%R_{\text{GGS}})_{21} \bigr]_{\alpha,\beta} = 8 \epsilon$, $\forall
%\alpha \prec \beta$. This result still holds if $J$ is altered by a
%different ordering in \eqref{jprod}.
%
(i) There exist unique half-integer $\kab$ such that \\ 
$\frac{d^2}{d\hbar^2} \biggl|_{\hbar=0} 
\bigl[\bar R_{J} - (\bar R_J)_{21}\bigr]_{\alpha,\beta} = 0, 
\forall \alpha \prec \beta.$

(ii) These $\kab$ are given by the combinatorial formula

\begin{multline} \label{fk}
K_{\alpha,\beta} = \frac{1}{2} [\alpha \lessdot \beta, \alpha \prec
\beta] - \frac{1}{2} [\alpha \gtrdot \beta, \alpha \prec \beta] +
[\exists \gamma \mid \alpha \prec \gamma \prec \beta, \alpha \lessdot
\gamma] \\ - [\exists \gamma \mid \alpha \prec \gamma \prec \beta,
\alpha \gtrdot \gamma] + [\alpha \prec^{\leftarrow}
\beta](1-|\alpha|).
\end{multline}

(iii)  For these $\kab$, and no others, one has $R_J \equiv R_{\text{GGS}}
\pmod{\hbar^3}$.
\end{thm}

{\it Proof.}  (i)  This is clear upon expanding $R_J$ modulo $\hbar^3$.  (See
Section 2 for details.)

(ii), (iii)  Proved in Section 2. \qs

\begin{conj} {\bf ``the twist conjecture''} \label{jconj} 
Taking $\kab$ as in \eqref{fk},

I. The matrix $R_J$ satisfies the QYBE.

II. The matrix $R_J$ coincides with $R_{\text{GGS}}$.
\end{conj}

The two parts of the twist conjecture are analogous to those of the GGS
conjecture in the following way.  Conjecture \ref{jconj}.II is a
strengthened version of \ref{ggs}.II, while Conjecture
\ref{jconj}.I is equivalent to \ref{ggs}.I modulo
\ref{jconj}.II.

\begin{thm} \label{comp} i) The twist conjecture holds for 
$n \leq 12$.  ii) The twist conjecture is true modulo $\hbar^3$.
\end{thm}

{\it Proof.}  (i) Given Theorem \ref{ggscomp}, it is sufficient to check
$J^{-1} \bar R_{\text{GGS}} J_{21} = \bar R_J$ for all triples, $n \leq 12$,
which has been carried out directly by computer.

(ii) Conjecture \ref{jconj}.I mod $\hbar^3$ is obvious from
construction, and Conjecture \ref{jconj}.II is true mod $\hbar^3$ as a
consequence of Theorem \ref{nkt}.\qs

\subsection{The generalized disjoint and Cremmer-Gervais triples} \label{it}

\begin{defe} \label{gdd} A triple $(\Gamma_1, \Gamma_2, \tau)$ is said to be
{\it generalized disjoint} if $\Gamma_1 = \bigcup_{i=1}^{m}
\Gamma_1^i$ where $\Gamma_1^i \perp \Gamma_1^j, i \neq j$ and
$\tau \Gamma_1^i \cap \Gamma_1 \subset \Gamma_1^{i+1}, i < m,$ and
$\tau \Gamma_1^m  \cap \Gamma_1 = \emptyset$.  If, in fact,
$\tau \Gamma_1^i \perp \Gamma_1^j, j \neq i+1$, and $\tau \Gamma_1^m \perp
\Gamma_1$, then the triple is said to be orthogonal generalized disjoint.
\end{defe}

\begin{exam} The case $\Gamma_1 = \{\alpha_i \mid i \not \equiv 0 \pmod 3,
i < n-3\}$, $\tau \alpha_i = \alpha_{i+3}$ is orthogonal generalized
disjoint.
\end{exam}

%In particular, the generalized disjoint case includes the {\it
%disjoint} case where $\Gamma_1 \cap \Gamma_2 = \emptyset$, discussed
%in \cite{H}.  

\begin{thm} \label{tdg}  (i) The twist conjecture is true in the disjoint case.
(ii) The twist conjecture is true in the orthogonal generalized disjoint
case.
\end{thm}

{\it Proof.}  See Sections \ref{ds}, \ref{gds}, and Appendices \ref{da}, 
\ref{gda}.  Note that the twist $\cal J$ used in
the disjoint case was first constructed by T. Hodges in \cite{H}.\qs

%Note that the twist $\cal J$ used in the disjoint case to prove $R_J$
%satisfies the QYBE was found by T. Hodges \cite{H}.

\begin{thm} The twist conjecture is true for the Cremmer-Gervais triples 
$(\Gamma_1, \Gamma_2, \tau)$ and $(\Gamma_2,\Gamma_1,\tau^{-1})$ where
$\Gamma_1 = \{\alpha_1, \ldots,\alpha_{n-2}\}, \Gamma_2 =
\{\alpha_2,\ldots,\alpha_{n-1}\}, \tau \alpha_i =
\alpha_{i+1}$.
\end{thm}

{\it Proof.}  In Section \ref{crgs}, we prove $R_J = R_{\text{GGS}}$.
On the other hand, it is known that $R_{\text{GGS}}$ satisfies the
QYBE in this case \cite{H2}.\qs

\begin{rem}
In fact, one may check $R_J=R_{\text{GGS}}$ when $\tau$ is replaced by
the map $\tau^k$ for $k \in \Z \setminus \{0\}$.  Combining this with
the result on generalized disjoint triples and a generalization of the
union arguments in the following section, one may conclude that $R_J =
R_{\text{GGS}}$ whenever $\Gamma_1 = \bigcup \Gamma_1^{(i)}$ where
$\Gamma_1^{(i)} \perp \Gamma_1^{(j)}, i \neq j$, and $\tau
\Gamma_1^{(i)} \cap \Gamma_1^{(j)} = \emptyset$ whenever $i > j$.
In particular, this includes the case when $\tau$ sends everything in
the same direction--i.e., $\tau(\alpha_i) = \alpha_j$ implies $j > i$ for all
$i$  (or $i < j$ for all $i$).  (The proof is omitted).
\end{rem}

\subsection{Maximal triples and unions}
In this subsection, we summarize reductions of the twist and
GGS conjectures which are proved in the following two subsections.

\begin{defe} \label{md} \cite{GGS} We say that $(\Gamma_1',\Gamma_2',\tau') < (\Gamma_1,\Gamma_2,\tau)$ if $\Gamma_1' \subset \Gamma_1$ and $\tau' = \tau \bigl|_{\Gamma_1'}$.
\end{defe}

The following theorem reduces the twist and GGS conjectures to the case of
maximal triples. 

\begin{thm} \label{mt}
Suppose $(\Gamma_1',\Gamma_2',\tau') < (\Gamma_1,\Gamma_2,\tau)$
are Belavin-Drinfeld triples.  Then if the twist or GGS conjecture holds
for the larger triple, it also holds for the smaller one.
\end{thm}

{\it Proof.}  See Section \ref{ms}. \qs

\begin{defe} \label{ud} Define $(\Gamma_1, \Gamma_2, \tau) = \bigcup (\Gamma_1^{(i)}, \Gamma_2^{(i)}, \tau^{(i)})$
by $\Gamma_1 = \bigcup \Gamma_1^{(i)}, \Gamma_2 = \bigcup
\Gamma_2^{(i)}$, with $\tau: \Gamma_1 \rightarrow \Gamma_2$ given by
$\tau\bigl|_{\Gamma_1^{(i)}} = \tau^{(i)}$.  Call a union {\it
orthogonal} if $\Gamma_1^{(i)} \perp \Gamma_1^{(j)}$ and
$\Gamma_2^{(i)} \perp \Gamma_2^{(j)}$.  Furthermore, an orthogonal
union is termed {\it $\tau$-orthogonal} if, in addition,
$\Gamma_2^{(i)} \cap \Gamma_1 \subset \Gamma_1^{(i)}, \forall i$.
\end{defe}

It is easy to check that an orthogonal union of Belavin-Drinfeld
triples always defines a Belavin-Drinfeld triple (one must check
that $\tau$ is nilpotent and a graph isomorphism for the union).
A triple that is a $\tau$-orthogonal union of two nonempty triples
is called {\it decomposable}; otherwise it is {\it indecomposable}.
The following theorem reduces the twist and GGS conjectures to the
case of indecomposable triples.

\begin{thm} \label{ut} If $(\Gamma_1, \Gamma_2, \tau)$ is a $\tau$-orthogonal
union of $(\Gamma_1^{(i)}, \Gamma_2^{(i)}, \tau^{(i)})$, then the twist or GGS
conjecture holds for the union iff it holds for each 
triple $(\Gamma_1^{(i)}, \Gamma_2^{(i)},
\tau^{(i)})$.
\end{thm}

{\it Proof.} See Section \ref{us}. \qs

\subsection{Maximal triples} \label{ms}

In this section, by Giaquinto's suggestion, we investigate the notion
of maximal triples as defined in Definition \ref{md} with the goal of
proving Theorem \ref{mt}.
We will assume throughout this section that $(\Gamma_1', \Gamma_2', \tau')
< (\Gamma_1, \Gamma_2, \tau)$.
Define $G, H, X$ for $(\Gamma_1, \Gamma_2, \tau)$ as in Section 2.1,
2.2 and similarly $G', H', X'$ for $(\Gamma_1', \Gamma_2', \tau')$.  
We begin with an important result.

\begin{prop} \label{mp1} If $\sum_{k=1}^m \bigl( \tau \alpha_{i_k} - \alpha_{i_k} \bigr) = 
\sum_{k=1}^m \bigl( \tau \alpha_{j_k} - \alpha_{j_k} \bigr)$ then
$\exists \rho \in S_m$, a permutation, such that $i_k = j_{\rho(k)},
\forall k$.
\end{prop}

{\it Proof.}  Suppose that $\not \exists \rho \in S_m$ such that $i_k
= j_{\rho(k)}, \forall k$.  Then, let $\alpha_l \in \Gamma_1$ be a
maximal simple root under the ordering $\prec$ so that $\tau \alpha_l
- \alpha_l$ does not appear the same number of times in the sequences
$(i_k), (j_k)$.  Then, $0 = \sum_{k=1}^m \bigl( \tau \alpha_{j_k} -
\tau \alpha_{i_k} - \alpha_{j_k} + \alpha_{i_k} \bigr) = p \tau
\alpha_l + \sum_k \alpha_{o_k}$ where $p \neq 0$ and $\tau \alpha_l
\neq \alpha_{o_k}$ for any $k$.  This is a contradiction.$\quad
\square$

Define $H = \{\sum_{i=1}^{m} \tau \alpha_{k_i} - \alpha_{k_i} \mid
\alpha_{k_i} \in \Gamma_1, m \geq 2 \}$, and $G = \{ e_{\alpha} \o
e_{\beta} \mid \alpha+\beta \in H$.  Clearly $(R_J)_+, (R_J)_- \in
\text{Span}_{\C[q,q^{-1}]} (G)$, and moreover, $G \subset \tilde
\Gamma_1 \times (\tilde \Gamma_2)^T \cup (\tilde \Gamma_2)^T \times
\tilde \Gamma_1$ where $T$ takes the transpose of any element.
Furthermore, let $V_0 \subset Mat_n(\C) \o Mat_n(\C)$ be the space of
zero weight in the representation $\g \o \g$ of $\g$.  Elements will
be said to have zero weight.  Then, we may define the following:

\begin{defe} For any element $x \in H$, define $i_\tau(x) = m$ if $x = 
\sum_{k=1}^m \bigl( \tau \alpha_{i_k} - \alpha_{i_k} \bigr)$. $m$ is 
called the {\it $\tau$-index} of $x$.  Similarly define 
$i_\tau (e_{\alpha} \o e_{\beta}) \equiv i_\tau(\alpha+\beta)$ 
when $e_{\alpha} \o e_{\beta} \in G$. 
%Further define
%$s_\tau(x) = \prod_{k=1}^m \bigl( \text{sign}(\alpha_{i_k},\tau
%\alpha_{i_k}) \bigr)$ and $s_\tau(e_\alpha \o
%e_\beta)=s_\tau(\alpha+\beta), \alpha+\beta \in H$.  $s_\tau(x)$ is
%called the {\it sign of $x$ with respect to $\tau$}.
\end{defe}

Clearly $i_\tau(x+y) = i_\tau(x)+ i_\tau(y)$ 
%s_\tau(x+y) = s_\tau(x)
%s_\tau(y)$ 
for $x,y \in H$ and $i_\tau(xy) = i_\tau(x) + i_\tau(y)$
%s_\tau(xy) = s_\tau(x) s_\tau(y)$ 
for $x,y \in G$ and $xy \neq 0$.
Note that if $|\tau^k \alpha|=|\alpha|$, then 
%$s_\tau(\tau^k \alpha-\alpha)=\text{sign}(\alpha,\tau^k \alpha)$ and 
$i_\tau(\tau^k \alpha-\alpha) = k|\alpha|$.  This concept of index
over $H$ and $G$ will come in handy in section 5.

\begin{prop} \label{mp} 
Suppose $(\alpha, \beta) \in X \setminus X'$ and $y \in G$.  Then
$\{\eab, e_{-\alpha} \otimes e_{\beta}, (\eab)y,$ \\ 
$(e_{-\alpha} \otimes e_{\beta})y, y(\eab), y(e_{-\alpha} \otimes
e_{\beta}\} \cap G' = \emptyset$.
\end{prop}

{\it Proof.}  Since $(\alpha, \beta) \notin X'$, it follows that
$\beta-\alpha \notin H'$, and hence $\beta-\alpha + H \cap H' = \emptyset$.
\qs

Denote by $K_{\alpha,\beta}$ and $K'_{\alpha,\beta}$ the appropriate
$K$-coefficients for the two triples, and by $\epsilon$ and
$\epsilon'$ the appropriate $\epsilon$-matrices.

\begin{corr} \label{mm} (i) If $(\alpha, \beta) \in X'$, then 
$K_{\alpha,\beta} = K'_{\alpha,\beta}$. 

Now suppose $e_{\alpha} \o e_{\beta} \in G'$.  Then (ii) $\epsilon_{\alpha, \beta} = \epsilon'_{\alpha, \beta}$, hence

(iii) $(R_J)_{\alpha, \beta} = (R'_J)_{\alpha, \beta}, \quad$ (iv) $(R_J)_{-\beta,-\alpha} = (R'_J)_{-\beta, -\alpha}$,

(v) $(R_{\text{GGS}})_{\alpha, \beta} = (R'_{\text{GGS}})_{\alpha, \beta}, 
\quad$ and (vi) $(R_{\text{GGS}})_{-\beta, -\alpha} = (R'_{\text{GGS}})_{-\beta, -\alpha}.$
\end{corr}

{\it Proof.} The proposition shows that terms from $G \setminus G'$
do not affect terms in $G'$ when expanding \eqref{nk}
and the matrices $\epsilon, R_J, R_{\text{GGS}}$.  Note that this result
also follows from the combinatorial formulas \eqref{fe}, \eqref{fk} (which
are derived independently). \qs

{\it Proof of Theorem \ref{mt}.}  (i) Choose $r^0$ to satisfy equation
\eqref{tr02} for the triple $(\Gamma_1,\Gamma_2,\tau)$.
Clearly the equations for $(\Gamma_1',\Gamma_2',\tau)$ are a subset of
these, and by Remark \ref{or0}, it is sufficient to consider only this
$r^0$.  With this $r^0$, define $\bar R_J, \bar
R_{\text{GGS}}$ corresponding to $(\Gamma_1, \Gamma_2, \tau)$ and
$\bar R'_J, \bar R'_{\text{GGS}}$ corresponding to $(\Gamma'_1,
\Gamma'_2, \tau)$.

Construct $f \in
\mathfrak{h}$ as follows:
\begin{align}
\tau \alpha_k f & = \alpha_k f, \quad \alpha_k \in \Gamma_1 \setminus \Gamma_1', \\
\tau \alpha_k f & = 1 + \alpha_k  f, \alpha_k \in \Gamma_1'.
\end{align}

Then, we see that $e^{ft} \otimes e^{ft} e_{\tau \alpha} \otimes
e_{-\alpha} e^{-ft} \otimes e^{-ft} = e^t e_{\tau \alpha} \otimes
e_{-\alpha}$ whenever $\alpha \in \Gamma_1 \setminus \Gamma_1'$, and
$e_{\tau \alpha} \otimes e_{-\alpha}$ otherwise (that is, when $\alpha
\in \tilde \Gamma_1'$).  This obviously holds as well for $e_{-\alpha}
\otimes e_{\tau \alpha}$. Clearly, if one takes any term $x$ of zero weight,
congugation by $e^{ft} \otimes
e^{ft}$ leaves the term unaltered.  This implies that $(e^{ft} \otimes
e^{ft}) x (e^{-ft} \otimes e^{-ft}) = e^{kt} x$ for $k \in Z^+$
whenever $x \in G \setminus G'$, but $k=0$ when 
$x \in G' \cup V_0$. It follows from Proposition \ref{mp} that
$\underset{t \rightarrow
-\infty}{\text{lim}} e^{ft} \otimes e^{ft} \bar R_J e^{-ft} \otimes
e^{-ft} = \bar R_J'$ and \\ $\underset{t \rightarrow -\infty}{\text{lim}}
e^{ft} \otimes e^{ft} \bar R_{\text{GGS}} e^{-ft} \otimes e^{-ft} =
\bar R'_{\text{GGS}}$.  This clearly implies the theorem. \qs

\subsection{Unions of triples} \label{us}
In this section, we investigate unions as defined in Definition \ref{ud} 
with the goal of proving
Theorem \ref{ut}. We will see in Lemma \ref{ul2}, in fact, that the
matrices $R_{\text{GGS}}$ and $R_J$ for the $\tau$-orthogonal union
follows directly from those for each triple.

Since it is clear that a union of triples is larger than each piece under
the ordering of the previous sections, we may pick $r^0$ to satisfy
\eqref{tr02} for the union which includes all equations for each
smaller triple.  Fix some such $r^0$, which will be sufficient by
Remark \ref{or0} to make statements about the conjectures.  We use the
notation $R^{(i)}_{\text{GGS}}$ and $R^{(i)}_J$ for the respective $R$
matrices.  Set $R'_s \equiv q^{\tilde r^0} R_s q^{\tilde r^0}$.  In
addition, set $S^a_b \equiv R^a_b - R'_s$ for any such subscripts $b$
and superscripts $a$, and similarly $\bar S^a_b \equiv \bar R^a_b -
R_s$.  We will use $H, G$ as defined in Section 2.2, and define
$H^{(i)}, G^{(i)}$ for the respective subtriples.  Finally, define
$V_k^{(i)} = \{e_{\alpha} \mid \alpha \in \tilde \Gamma_k^{(i)} \cup
\tilde \Gamma_2^{(i)} \}$ for $k \in \{1,2\}$, and let $V^{(i)} =
\{e_{\alpha} \mid \alpha \in (\tilde \Gamma_1 + \tilde \Gamma_2) \cap
\Gamma \}$.

\begin{lemma} \label{ul} Suppose that $(\Gamma_1, \Gamma_2, \tau)$ is
a $\tau$-orthogonal union of $(\Gamma_1^{(i)},
\Gamma_2^{(i)})$.  Then $V_k^{(i)} V_k^{(j)} = \emptyset$ for $i \neq j,
k \in \{1,2\}$, and $V_1^{(i)} (V_2^{(j)})^T = (V_2^{(j)})^T V_1^{(i)} =
\emptyset$ for $i \neq j$, where $T$ takes the transpose of each element.
Hence one has $(G^{(i)})_{ab} X (G^{(j)})_{cd} = \emptyset$ for $a,b,c,d \in
\{1,2,3\}, a \neq b, c \neq d$ and any $X \in V_0 \o V_0 \o V_0$.
\end{lemma}

{\it Proof.} Indeed, the first two assertions follow from the facts
$\Gamma_1^{(i)} \perp \Gamma_1^{(j)}, i \neq j$, and $\Gamma_1^{(i)} \cap
\Gamma_2^{(j)} = \emptyset, i \neq j$.  The remainder follows since $G^a_b V_0, V_0 G^a_b \subset \text{Span}(G^a_b)$.
$\quad \square$

\begin{lemma} \label{ul2} Suppose that $(\Gamma_1, \Gamma_2, \tau)$ is
a $\tau$-orthogonal union of $(\Gamma_1^{(i)}, \Gamma_2^{(i)})$.
Then $S_t = \sum_i S^{(i)}_t$ for $t \in \{\text{GGS}, J\}$.
\end{lemma}

{\it Proof.}
Clearly any element of $G$ is of the form $(\eab)_{ab}$
for $\{a,b\} = \{1,2\}$ and $\alpha \in \tilde \Gamma_1, \beta \in
\tilde \Gamma_2$.  It is clear that $\tilde \Gamma_1 = \sqcup_i \tilde
\Gamma_1^{(i)}$ where $\sqcup$ denotes a disjoint union. Since we also
have $\tilde \Gamma_2^{(i)} \cap \tilde \Gamma_1^{(j)} = 0$ for $i
\neq j$, it follows that $G = \sqcup_i G^{(i)}$. Hence, $S_t =
\sum_i S^{(i)}_t$ for $t \in \{\text{GGS}, J\}$. \qs

\begin{thm} If $(\Gamma_1, \Gamma_2, \tau)$ is a $\tau$-orthogonal
union of $(\Gamma_1^{(i)}, \Gamma_2^{(i)})$, then any part of the twist
or GGS conjecture holds for each triple $(\Gamma_1^{(i)},
\Gamma_2^{(i)}, \tau^{(i)})$ iff the same part of the twist or GGS
conjecture holds for the union.
\end{thm}

{\it Proof.}  By Lemma \ref{ul}, we have $(S_t^{(i)})_{ab} X
(S_t^{(j)})_{cd} = 0$ for $a,b,c,d \in \{1,2,3\}, a \neq b, c \neq d,$
and any $X \in \text{End}(Mat_n(\C) \o Mat_n(\C) \o Mat_n(\C))$, one
may consider the QYBE separately in $V_0 \o V_0 \o V_0$ and in
$V^{(i)} \o V^{(i)} \o V^{(i)}$ for each $i$, and one may consider the
Hecke relation and $R_J = R_{\text{GGS}}$ separately in $V_0 \o V_0,
V^{(i)} \o V^{(i)}$.  The $V_0$ components clearly only involve $R_s$
so are satisfied, while each $V^{(i)}$ component holds iff the
respective equation holds for the $i$-th subtriple.  Finally, the
respective equation holds for the union iff it holds in each
component, and each component yields the same equation considered in
the union and in the appropriate subtriple by Lemma \ref{ul2}.  The
theorem follows from these observations. \qs

\section{Complete description of $K_{\alpha,\beta}$}

In this section we prove Proposition \ref{fep} and Theorem \ref{nkt} by first
deducing \eqref{fk} from $\frac{d^2}{d\hbar^2} \bigl|_{\hbar=0} \bar
R_{\text{GGS}}$ and then applying the development to prove the
equivalence of \eqref{eps} and \eqref{fe} for $\epsilon$.  To do this,
we will rely on important bijections between different ways the product
of two terms in $aP,Pa,a^2$ can arise, where corresponding pairs
cancel in the expansion.

%is the unique $K_{\alpha,\beta}$ such that $\bar R_J
%\equiv \bar R_{\text{GGS}}$
%according to simple combinatorial information about pairs $\alpha \prec
%\beta$ that appear in Belavin-Drinfeld triples.

\begin{lemma} \label{ch3} There exist unique $\kab$ such that $\frac{d^2}{d\hbar^2}
(R_J - (R_J)_{21})_{\alpha,\beta} = 0, \alpha \prec \beta$.  These are
given by the formula
\begin{equation} \label{nk}
K_{\alpha,\beta} = \text{sign}(\alpha,\beta) \biggl[ \sum_{i \geq
j} a_+^{i} a_+^{j} - \sum_{i < j} a_+^{i} a_+^{j} + a_+ P_- +
a_- P_+ + P_+ a_+ + a_+ a_- - a_- a_+ \biggr]_{\alpha,\beta}.
\end{equation}
With these $K_{\alpha,\beta}$, the condition that $R_J \equiv R_{\text{GGS}}
\pmod{\hbar^3}$ reduces to showing the equivalence of \eqref{nk}
and \eqref{fk}, and that $\kab = \epsilon_{\alpha,\beta} = 0$
when $\alpha \not \prec \beta$, in this case defining $\kab$ by \eqref{nk} for all
$\alpha,\beta \in \tilde \Gamma$.
\end{lemma}

{\it Proof.}  We expand $\bar R_J$ modulo $\hbar^3$
as follows:
\begin{multline} \label{rjh3}
\bar R_J \equiv \biggl[ \sideset{}{^>}\prod_{\pab \in X} \bigl( 1 -
2\: \text{sign}(\alpha, \beta) \hbar(1 + K_{\alpha,\beta} \hbar)\:
e_\beta \otimes e_{-\alpha} \bigr) \biggr] \biggl[ 1 +
\frac{\hbar^2}{2} \sum_{i} e_{ii} \otimes e_{ii} + 2 \hbar P_- \biggr]
\\ \biggl [ \sideset{}{^<}\prod_{\pab \in X} \bigl( 1 + 2\:
\text{sign}(\alpha, \beta)\hbar(1 + K_{\alpha,\beta} \hbar)\:
e_{-\alpha} \otimes e_{\beta} \bigr) \biggr] \equiv 1 + 2\hbar \biggl
[ P_- + \sum_{\alpha \prec \beta} \text{sign}(\alpha, \beta)\:
e_{-\alpha} \wedge e_\beta \biggr] \\ + \hbar^2 \biggl[ \frac{1}{2}
\sum_i e_{ii} \otimes e_{ii} + 2 \sum_{\alpha \prec \beta} K_{\alpha,
\beta}\: \text{sign}(\alpha, \beta)\: e_{-\alpha} \wedge e_{\beta} + 4
a_+ P_- + 4 P_- a_- + 4 a_+ a_- \\ + 4 \sum_{i \geq j} a_+^{i} a_+^{j}
+ 4 \sum_{i < j} a_-^{i} a_-^{j} \biggr] \pmod {\hbar^3}.
\end{multline}

\noindent
If we skew-symmetrize the second order terms in \eqref{rjh3} by the
substitution $x \mapsto x - x_{21}$, \eqref{nk} follows
as a necessary and sufficient condition for
$\frac{d^2}{\hbar^2}\bigl|_{\hbar = 0} \bigl[ \bar R_{J} - (\bar
R_{J})_{21}\bigr]_{\alpha,\beta} = 0.$

It is clear that $R_J \equiv R_{\text{GGS}} \equiv 1+2\hbar r \pmod{\hbar^3}$,
so in particular by the comments in Section 1.2, this implies $\frac{1}{2}(R_J + (R_J)_{21}) \equiv R_{\text{GGS}} \equiv 1 + \hbar P + 2 \hbar^2 s \pmod{\hbar^3}$.  All that remains, then, is to show that \eqref{fk} and \eqref{nk} are
equivalent, and that $\epsilon_{\alpha,\beta} = \epsilon_{-\beta,-\alpha} = 0$
when $\alpha \not \prec \beta$ and the same for $K_{\alpha,\beta}$ in 
\eqref{nk}.\qs

Take positive roots $\alpha, \beta$, $|\alpha| = |\beta|, \alpha \neq
\beta$.  Set $\alpha = e_i - e_j$, $\beta = e_k - e_l$.  Then, say
that $\alpha < \beta$ if $i < k$, and in this case, $\alpha \lessdot
\beta$ if $j = k$, $\alpha\: \overline{<}\: \beta$ if $j > k$, and
$\alpha \ll \beta$ if $j < k$ (we repeat the definition of $\lessdot$
given in Section 1 for completeness.)  We will use $\alpha > \beta$ if
$\beta < \alpha$, and similarly, $\gtrdot, \overline{>} ,$ and $\gg$
are the reverse directions of $\lessdot, \overline{<},$ and $\ll$,
respectively.  With these definitions, we will take $x^{<} =
\sum_{\alpha < \beta} \bigl(x_{\alpha,\beta}\: \eab +
x_{-\beta,-\alpha}\: e_{-\alpha} \otimes e_{\beta} \bigr)$ and
similarly for the other defined relations.

\begin{lemma} \label{bkl} We may rewrite \eqref{nk} as follows:
\begin{multline} \label{bk}
K_{\alpha,\beta} = \text{sign}(\alpha,\beta) \biggl[ \frac{1}{2}
\bigl(a_+^{\gtrdot} - a_+^{\lessdot} \bigr) + a_+^{\overline{>}} P - P
a_+^{\overline{<}} + \sum_{\alpha' \prec^\leftarrow \beta'}
\text{sign}(\alpha', \beta') (1 - |\alpha'|)\: e_{\beta'} \otimes
e_{-\alpha'} \\ + \sum_{i < j} \bigl( a_+^{j} a_+^{i} - a_+^{i}
a_+^{j} \bigr) +a_+ a_- - a_- a_+ \biggr]_{\alpha, \beta}.
\end{multline}
\end{lemma}

{\it Proof.} First, we expand $a_+ P_- + a_- P_+ + P_+ a_+$:
\begin{multline} \label{ap}
\bigl( a_+ P_- + a_- P_+ + P_+ a_+ \bigr)_+ = \biggl[a_+^{\overline{>}} P +
\frac{1}{2} a_+^{\gtrdot} + a_-^{<} P + P a_+^{\ll} + 
\frac{1}{2} a_+^{\lessdot}\biggr]_+ = \\
\biggr[ \frac{1}{2} \bigl(
a_+^{\lessdot} + a_+^{\gtrdot} \bigr) + a_+^{\overline{>}} P - P 
a_+^{<} + P a_+^{\ll} \biggr]_+ =
\frac{1}{2} \bigl(a_+^{\gtrdot} - a_+^{\lessdot} \bigr) + 
a_+^{\overline{>}} P - P a_+^{\overline{<}}.
\end{multline}

Now we simplify $\sum_{i} (a_+^{i})^2$.  $\sum_{i} (a_+^{i})^2 =
\sum_{\alpha \prec \beta, \alpha' \prec \beta'} \text{sign}(\alpha,
\beta) \text{sign}(\alpha', \beta')\: e_{\beta} e_{\beta'} \otimes
e_{-\alpha} e_{-\alpha'}$. For $\alpha, \alpha' \in \tilde \Gamma_1$,
$\tau^i(\alpha) = \beta, \tau^i(\alpha') = \beta'$, one sees that
$e_{\beta} e_{\beta'} = e_{\beta + \beta'}$ and $e_{-\alpha}
e_{-\alpha'} = e_{-\alpha - \alpha'}$ iff \\ $\alpha+\alpha' \in
\tilde \Gamma_1$ and $\tau^i(\alpha+\alpha') = \beta+\beta'$, reversing
order.  Thus, since in this case \\ $\text{sign}(\alpha, \beta)
\text{sign}(\alpha', \beta') = -\text{sign}(\alpha + \alpha', \beta +
\beta')$,

\begin{equation} \label{a2}
\sum_{i} (a_+^{i})^2 = -\sum_{\alpha \prec^\leftarrow \beta} \text{sign}(\alpha, \beta) (|\alpha| - 1)\: e_{\beta} \otimes e_{-\alpha}.
\end{equation}

It is clear that \eqref{ap} and \eqref{a2} imply the proposition.
$\quad \square$

Now, we proceed to show the equivalence of \eqref{bk} and \eqref{fk}
by canceling most terms in the expansion of \eqref{nk} pairwise.  Define
the following sets:

\begin{gather}
M_1 = \{ ((\alpha, \tau^x \alpha), (\beta, \tau^y \beta)) \in X
\times X \mid \alpha \gtrdot \beta, \tau^x \alpha \lessdot \tau^y
\beta, x > y\}, \\
M_2 = \{ ((\alpha, \tau^x \alpha), (\beta, \tau^y
\beta)) \in X \times X \mid \alpha \gtrdot \beta, \tau^x \alpha
\lessdot \tau^y \beta, x < y\}, \\
M_3 = \{ ((e_x - e_y, e_u-e_v),
(e_{v'} - e_v, e_x - e_{x'})) \in X \times X \mid x' < y, u < v'\}, \\ 
M_4 = \{ ((e_x - e_y, e_u-e_v), (e_u - e_{u'}, e_{y'} - e_y)) \in
X \times X \mid x < y', u' < v\}, \\
M_5 = \{(\alpha,\beta) \in X \mid \alpha \overline{>} \beta\}, \\ 
M_6 = \{(\alpha,\beta) \in X \mid \alpha \overline{<} \beta\}.
\end{gather}

Clearly these are defined so that the following hold:

\begin{gather}
\sum_{i < j} a_+^j a_+^i =
\sum_{((\alpha,\beta), (\gamma,\delta)) \in M_1} a_{\alpha,\beta}
a_{\gamma,\delta} e_{\beta+\delta} \o e_{-\alpha-\delta}, \\
\sum_{i < j} a_+^i a_+^j =
\sum_{((\alpha,\beta), (\gamma,\delta)) \in M_2} a_{\alpha,\beta}
a_{\gamma,\delta} e_{\beta+\delta} \o e_{-\alpha-\delta}, \\
(a_+ a_-)_+ = \sum_{((\alpha,\beta),(\gamma,\delta)) \in M_3} a_{\alpha,\beta}
a_{-\delta,-\gamma} e_{\alpha-\delta} \o e_{\beta-\gamma}, 
\end{gather}
\begin{gather}
(a_- a_+)_+ = \sum_{((\alpha,\beta),(\gamma,\delta)) \in M_4} a_{-\delta,-\gamma}
a_{\alpha,\beta} e_{\alpha-\delta} \o e_{\beta-\gamma}, \\
a_+^{\overline{>}} P = \sum_{(\alpha,\beta) \in M_5} a_{\alpha,\beta} \eab P, \quad
P a_+^{\overline{<}} = \sum_{\alpha,\beta \in M_6} a_{\alpha,\beta} P \eab.
\end{gather}

%$ and
%similarly for $M_2$ and $\sum_{i<j} a_+^i a_+^j$.  Then, we have
%$a_+ a_- = \sum_{((\alpha,\beta),(\gamma,\delta))} a_{\alpha,\beta}
%a_{-\delta,-\gamma}$, $
%$M_3$ and $a_+ a_-$; and $M_4$ and $a_- a_+$.  Furthermore, it is easy
%to see $M_5$ corresponds to terms in $a^{\overline{>}} P$, and $M_6$
%corresponds to those in $P a^{\overline{<}}$.  
Now we define subsets
$M_i' \subset M_i$, set $M_i'' = M_i \setminus M_i'$, and bijections
$f: M_1' \rightarrow M_3', g: M_2' \rightarrow M_4', f': M_1''
\rightarrow M_5', g': M_2'' \rightarrow M_6'$ which allow pairwise
cancellation, leaving us to expand \eqref{bk} by only those terms in
$M_3'',M_4'',M_5'',M_6''$, which will lead directly to \eqref{fk}.
Define $M_i'$ as follows:

\begin{gather}
M_1' = \{ ((\alpha, \tau^x \alpha), (\beta, \tau^y \beta)) \in M_1
\mid x-y \nmid y \}, \\
M_2' = \{ ((\alpha, \tau^x \alpha), (\beta,
\tau^y \beta)) \in M_2 \mid y-x \nmid x \}, \\
M_3' = \{ ((e_x -
e_y, e_u-e_v), (e_{v'} - e_v, e_x - e_{x'})) \in M_3 \mid (y-x') 
\nmid (y-x) \}, \\
M_4' = \{ ((e_x - e_y, e_u-e_v), (e_u - e_{u'},
e_{y'} - e_y)) \in M_4 \mid (y'-x) \nmid (y-x)\}, \\
M_5' = \{(e_x - e_y,e_u - e_v) \in M_5 \mid
x-u \nmid y-x \}, \\
M_6' = \{(e_x - e_y, e_u - e_v) \in M_6 \mid u-x \nmid y-x \}.
\end{gather}

Now, we construct bijections $f,g,f',g'$.  We begin with $f$.  Take
$((\alpha, \tau^x \alpha), (\beta, \tau^y \beta)) \in M_1'.$ Suppose
$y = p(x-y) + q$ where $p,q \in \N$ and $0 < q < x-y$.  Then $\alpha
\gtrdot \beta \gtrdot \tau^{x-y}(\alpha+\beta) \gtrdot \cdots \gtrdot
\tau^{p(x-y)}(\alpha+\beta) \gtrdot \tau^{(p+1)(x-y)} \alpha$.  Then,

\begin{multline} \label{fdef}
f((\alpha, \tau^x \alpha),
(\beta, \tau^y \beta)) = \bigl[\bigl(
(1+\tau^{x-y}+\ldots+\tau^{p(x-y)})(\alpha+\beta)+\tau^{(p+1)(x-y)} \alpha,
\\
(\tau^q+\tau^{q+(x-y)}+\ldots+\tau^y)(\alpha+\beta)+\tau^x \alpha \bigr),
\bigl(
(\tau^{q}+\ldots+\tau^{q+(p-1)(x-y)})(\alpha+\beta)+\tau^y \alpha,
\\
(\tau^{x-y} + \ldots + \tau^{p(x-y)})(\alpha+\beta)+\tau^{(p+1)(x-y)} \alpha
\bigr)\bigr] \in M_3'.
\end{multline}

Similarly, if $((\alpha, \tau^x \alpha), (\beta, \tau^y \beta)) \in
M_2',$ one sets $x = p(y-x) + q$, $0 < q < x-y$, notices $\beta
\lessdot \alpha \lessdot \tau^{y-x}(\alpha+\beta) \lessdot \cdots
\lessdot \tau^{p(y-x)}(\alpha+\beta) \lessdot \tau^{(p+1)(x-y)}
\beta$, and is able to define

\begin{multline} \label{gdef}
g((\alpha, \tau^x \alpha),
(\beta, \tau^y \beta)) = \bigl[\bigl(
(1+\tau^{y-x}+\ldots+\tau^{p(y-x)})(\alpha+\beta)+\tau^{(p+1)(y-x)} \beta,
\\
(\tau^q+\tau^{q+(y-x)}+\ldots+\tau^x)(\alpha+\beta)+\tau^y \beta \bigr),
\bigl(
(\tau^{q}+\ldots+\tau^{q+(p-1)(y-x)})(\alpha+\beta)+\tau^x \beta,
\\
(\tau^{y-x} + \ldots + \tau^{p(y-x)})(\alpha+\beta)+\tau^{(p+1)(y-x)} \beta
\bigr)\bigr] \in M_4'.
\end{multline}

Next, we define $f'$ and $g'$:

\begin{gather}
f'((e_j - e_i,e_a - e_{a+i-j}),(e_k - e_j,e_{a+i-j} - e_{a+i-k})) = 
(e_{a+i-k} - e_i,e_a - e_k) \in M_5'', \\
g'((e_j - e_i, e_a - e_{a+i-j}),(e_k-e_j, e_{a+i-j}-e_{a+i-k})) = (e_k - e_a, e_i - e_{a+i-k}) \in M_6''.
\end{gather}

\begin{lemma} \label{mcl} (i) $f: M_1' \rightarrow M_3'$ is bijective. Given any
$((\alpha,\beta),(\gamma,\delta)) \in M_1'$, and \\
$f((\alpha,\beta),(\gamma,\delta)) =
((\alpha',\beta'),(\gamma',\delta')) \in M_3'$, one has
$a_{\alpha,\beta} a_{\gamma,\delta} + a_{\alpha',\beta'}
a_{-\delta',-\gamma'} = 0$.

(ii) $g: M_2' \rightarrow M_4'$ is bijective. Given any
$((\alpha,\beta),(\gamma,\delta)) \in M_2'$, and
$g((\alpha,\beta),(\gamma,\delta)) =
((\alpha',\beta'),(\gamma',\delta')) \in M_4'$, one has
$a_{\alpha,\beta} a_{\gamma,\delta} + a_{\alpha',\beta'}
a_{-\delta',-\gamma'} = 0$.

(iii) $f': M_1'' \rightarrow M_5'$ is bijective. Given any
$((\alpha,\beta),(\gamma,\delta)) \in M_1''$, and
$f'((\alpha,\beta),(\gamma,\delta)) = (\alpha',\beta') \in M_5'$, one
has $a_{\alpha,\beta} a_{\gamma,\delta} + a_{\alpha',\beta'} = 0$.

(iv) $g': M_2'' \rightarrow M_6'$ is bijective. Given any
$((\alpha,\beta),(\gamma,\delta)) \in M_2''$, and
$g'((\alpha,\beta),(\gamma,\delta)) = (\alpha',\beta') \in M_5'$, one
has $a_{\alpha,\beta} a_{\gamma,\delta} + a_{\alpha',\beta'} = 0$.
\end{lemma}

{\it Proof.} (i) Take any $((e_x - e_y, e_u-e_v), (e_{v'} - e_v, e_x -
e_{x'})) \in M'_3$. We find its inverse under $f$ and verify the identity.  
Suppose $\tau^r (e_x-e_y) = e_u - e_v, \tau^s
(e_{v'}-e_v) = e_x - e_{x'}.$ Clearly $\tau^r$ preserves orientation
on $e_x-e_y$.

Suppose $\tau^s$ reverses orientation on $e_{v'}-e_v$.
In this case, nilpotency of $\tau$ shows that $x'-x \leq y-x'$,
so $y-x' \nmid x'-x$ implies $x'-x < y-x'$.  Then,
one sees that $\tau^{2r+s}(e_{y-(x'-x)}-e_{y}) = e_u-e_{u+(x'-x)}$,
while $\tau^r(e_x'-e_{y-(x'-x)}) = e_{u+(x'-x)}-e_{u+(y-x')}$.  It
is easy to check that $f((e_{y-(x'-x)}-e_y,e_u-e_{u+(x'-x)}),
(e_{x'-e_{y-(x'-x)}},e_{u+(x'-x)}-e_{u+(y-x')}))=((e_x-e_y,e_u-e_v),(e_{v'}-e_v, e_x - e_{x'}))$ as desired.  Furthermore, we see that $\tau^{2r+s}, \tau^s$
reverse orientation while $\tau^r$ preserves orientation, so the desired
identity follows.

Now, suppose $\tau^s$ preserves orientation on $e_{v'} - e_v$. Then, $\tau^{r+s}(\alpha_i)
= \alpha_{i-(y-x')}$ for $x+y-x' \leq i \leq y$.  Then, suppose
$x'-x = p(y-x') + q$, $0 < q < y-x'$.  In this case,
$\tau^{(p+1)(r+s)+r}(e_{y-q}-e_y)=(e_u-e_{u+q})$ and
$\tau^{p(r+s)+r}(e_{x'}-e_{y-q})=(e_{u+q}-e_{v'})$. One
may check $f((e_{y-q}-e_y,e_u-e_{u+q}),(e_{x'}-e_{y-q},
e_{u+q}-e_{v'})) = ((e_x -
e_y, e_u-e_v), (e_{v'} - e_v, e_x - e_{x'}))$, as desired.  Since
$\tau^{r}, \tau^s$ both preserve orientation, the desired identity follows.

(ii)  This follows exactly as in (i).

(iii) Take any $(e_x - e_y, e_u - e_v) \in M_5'$.  We find its inverse
under $f'$ and verify the identity.  Indeed, suppose $\tau^r(e_x-e_y)
= e_u-e_v$ and $y-x = p(x-u) + q$ for $0 < q < x-u$.  Then,
$\tau^{(p+1)r}(e_{y-q} - e_y) = e_u - e_{u+q}$, and $\tau^{pr}(e_v - e_{y-q})
= e_{u+q} - e_x$, so that $f'((e_{y-q}-e_y),(e_u-e_{u+q})) = ((e_v-e_{y-q}),
(e_{u+q} - e_x))$.  Furthermore, $\tau^r$ preserves orientation on $e_x - e_y$
so the identity is verified.

(iv)  This follows exactly as in (iii).\qs

\begin{prop} (i) Formula \eqref{bk} is equivalent to \eqref{fk}.
(ii) Formula \eqref{eps} for $\epsilon$ is equivalent to \eqref{fe}.
\end{prop}

{\it Proof.}  (i) Given any $\tau^z \alpha = \beta$, it is clear that
$\exists ((\gamma,\delta),(\gamma',\delta')) \in M_3''$ with
$\gamma-\delta' = \alpha, \delta-\gamma' = \beta$ iff $\exists t, 0 <
t < z, t \nmid z$ such that $\tau^t \alpha \lessdot \alpha$.  In this
case it is easy to see (along similar lines as (i) in the proof of
Lemma \ref{mcl}) that $a_{\gamma,\delta} a_{-\delta',-\gamma'} =
-\text{sign}(\gamma',\delta')= -\text{sign}(\alpha,\beta)$.
Similarly, given any $\tau^z \alpha = \beta$, $\exists ((\gamma,
\delta),(\gamma',\delta')) \in M_4''$ such that $\gamma'-\delta =
\alpha, \delta'-\gamma = \beta$ iff $\exists t, 0 < t < z, t \nmid z$
such that $\tau^t \alpha \gtrdot \alpha$.  In this case,
$a_{-\delta,-\gamma} a_{\gamma',\delta'} = -\text{sign}(\gamma,\delta)
= -\text{sign}(\alpha,\beta)$.  Next, we find that $M_5'' = \{ (\alpha,
\tau^{tk} \alpha) \mid t,k \in \Z^+, \tau^t \alpha \lessdot \beta \}$
and $M_6'' = \{ (\alpha, \tau^{tk} \alpha) \mid t,k \in \Z^+, \tau^t
\alpha \gtrdot \beta \}$, and it is clear that all terms in $a_+^{\overline{>}}, a_+^{\overline{<}}$ appear with coefficient $-1$.  Hence, combining
Lemma \ref{mcl} with \eqref{bk}, we obtain precisely \eqref{fk}. 

(ii) Indeed, using $r_s = P + P_- - P_+$ we find $\epsilon_+ =
a_+ a_+ + a_+ a_- + a_- a_+ + P a_+^{\overline{<}} +
a_+^{\overline{>}} P + \frac{1}{2} \bigl( a_+^{\lessdot} +
a_+^{\gtrdot} \bigr)$.  Then, we see as in Lemma \ref{bkl} that
$\sum_i (a_+^{i})^2 = \sum_{\alpha \prec^{\leftarrow} \beta}
(1-|\alpha|) \text{sign}(\alpha,\beta)$, so $(a_+
a_+)_{\alpha,\beta} = [\alpha \prec^{\leftarrow}
\beta](-1)^{1-|\alpha|}(1-|\alpha|) + \sum_{i < j} (a_+^i a_+^j +
a_+^j a_+^i)$.  Hence, \eqref{fe} follows from the observations in (i). \qs

Proposition \ref{fep} and Theorem \ref{nkt} are proved. \qs

\section{The disjoint case and its generalization}

\subsection{The disjoint case} \label{ds}
This section is devoted to proving the following theorem:

\begin{thm}\label{dt} If $\Gamma_1 \cap \Gamma_2 = \emptyset$, then 
$R_{\text{GGS}} = R_J$.
\end{thm}

{\it Proof.} We will assume $\Gamma_1 \cap \Gamma_2 = \emptyset$ throughout this section.
The first observation to make in this case is that, since $\tau^2 =
0$, $J = J^{1} = 1 + A^1 = 1 + (q - q^{-1})
\sum_{\alpha \prec \beta}
\text{sign}(\alpha, \beta) q^{K_{\alpha,\beta}}\: e_\beta \otimes
e_{-\alpha}$.    Set $A=A^1$.

Let $\alpha \perp \beta$ denote either $\alpha \ll \beta$ or $\beta
\ll \alpha$ (this is {\bf not} the same as $(\alpha,\beta) = 0$). Using
\eqref{fk}, the form of $K_{\alpha, \beta}$ in the disjoint case is
summarized in the following table:

\begin{center}
\begin{tabular}{|c||c|c|}\hline
$K_{\alpha,\beta}$ & $\alpha \prec^\rightarrow \beta$ & $\alpha
\prec^\leftarrow \beta$ \\ \hline \hline 
$\alpha \perp \beta$ & 0 & $1 - |\alpha|$ \\ \hline 
$\alpha 
\lessdot \beta$ & $\frac{1}{2}$ & $\frac{3}{2} - |\alpha|$ \\ \hline
$\alpha \gtrdot \beta$ & $-\frac{1}{2}$ & $\frac{1}{2} -
|\alpha|$ \\ \hline
\end{tabular} \\ \vskip 8 pt
Table 4.1: $K_{\alpha, \beta}$ in the disjoint case.
\end{center}

Also, $\epsilon$ is summarized as follows:

\begin{center}
\begin{tabular}{|c||c|c|}\hline
$\epsilon_{\alpha,\beta}$ & $\alpha \prec^\rightarrow \beta$ & $\alpha
\prec^\leftarrow \beta$ \\ \hline \hline 
$\alpha \perp \beta$ & 0 & $(-1)^{1-|\alpha|} (1-|\alpha|)$ \\ \hline 
$\alpha 
\lessdot \beta$ & $-\frac{1}{2}$ & $(-1)^{1-|\alpha|}(\frac{1}{2}-|\alpha|)$ \\ \hline
$\alpha \gtrdot \beta$ & $-\frac{1}{2}$ & $(-1)^{1-|\alpha|}(\frac{1}{2}-|\alpha|)$ \\ \hline
\end{tabular} \\ \vskip 8 pt
Table 4.2: $\epsilon$ in the disjoint case.
\end{center}

Set $B = J^{-1} - 1$. Then, the following lemma describes $\bar R_J$:

\begin{lemma} \label{brjd}
(i) $B$ is given by a sum $\sum_{\alpha \in \tilde \Gamma_1} B_{\alpha,\tau \alpha}\:
e_{\tau \alpha} \otimes e_{-\alpha}$.

(ii) $\bar R_J$ is given by the following equation:
\begin{equation} \label{jdds}
\bar R_J = R_s + B + A_{21} + (q-1)B^{\gtrdot} + (q^{-1}-1)A_{21}^{\lessdot}
\end{equation}
\end{lemma}

{\it Proof.} (i) Note that $e_{-\alpha} e_{-\beta} = e_{-\alpha-\beta}$, $e_{\tau \alpha} e_{\tau \beta} = e_{\tau \alpha + \tau \beta} = e_{\tau(\alpha + \beta)}$, for $\alpha, \beta \in \tilde \Gamma_1$. Thus, when $B$ is
expanded, all terms will remain of this type.

(ii) It is clear that $\bar R_J = (1 + B)R_s(1 + A_{21}) = R_s + B R_s
+ R_s A_{21} + B R_s A_{21}$.  Since $e_{\beta} e_{-\alpha} =
e_{-\alpha} e_{\beta} = 0$ for all $\beta \in \tilde \Gamma_2, \alpha \in
\tilde \Gamma_1$, we see that $B R_s A_{21} = (q - q^{-1}) B \bigl(
\sum_{\alpha > 0} e_{-\alpha} \otimes e_{\alpha} \bigr) A_{21}$.
Also, $B \bigl( \sum_{\alpha > 0} e_{-\alpha} \otimes e_{\alpha}
\bigr) = PB^{<}_{21}$ and $\bigl( \sum_{\alpha > 0} e_{-\alpha}
\otimes e_{\alpha} \bigr) A_{21} = P(A_{21}^< - A_{21}^\lessdot)$. So,

\begin{multline*}
\bar R_J = R_s + B R_s + R_s A_{21} + (q - q^{-1}) P B^{<}_{21}
A_{21} \\ = R_s + (q - q^{-1}) P\biggl( B^{<}_{21} + A_{21}^{<}
- A_{21}^{\lessdot} + B^{<}_{21} A_{21}^{<} \biggr) + B +
A_{21} + (q-1)(B^{\gtrdot} + A_{21}^{\lessdot}).
\end{multline*}

Since $B^{<}_{21} + A_{21}^{<} + B^{<}_{21} A_{21}^{<} =
(1 + A_{21}^{<})(1 + B^{<}_{21}) - 1 = 0$, 
we find:

\begin{multline*}
\bar R_J = R_s - (q - q^{-1}) P A_{21}^{\lessdot} + B +
A_{21} + (q-1)(B^{\gtrdot} + A_{21}^{
\lessdot}) \\ = R_s + B + A_{21} + (q-1)B^{\gtrdot}
+ (q^{-1}-1)A_{21}^{\lessdot}.
\end{multline*}

The lemma is proved.$\quad \square$

Now, we compute $B$ using \eqref{jiprod}, in which $(\alpha, \beta) > (\alpha', \beta')$ whenever \\ $(e_{\beta} \otimes e_{-\alpha})
(e_{\beta'} \otimes e_{-\alpha'}) \neq 0$.  Define $L_{\alpha,
\beta}$ as follows:
\begin{equation}
L_{\alpha, \beta} = 
\begin{cases}
0 & \text{if}\ \alpha \perp \beta, \\ 
\frac{1}{2} & \text{if}\ \alpha \lessdot \beta, \\ 
-\frac{1}{2} & \text{if}\ \alpha \gtrdot \beta.
\end{cases}
\end{equation}

\begin{lemma}
(i) If $\alpha \prec^{\rightarrow} \beta$, then $B_{\alpha, \beta} = -A_{\alpha, \beta} = -(q - q^{-1}) q^{K_{\alpha,\beta}}$.

(ii) If $\alpha \prec^{\leftarrow} \beta$, then $B_{\alpha, \beta} = -\text{sign}(\alpha, \beta) (q - q^{-1})
q^{|\alpha|-1+L_{\alpha,\beta}}$.
\end{lemma}

{\it Proof.} (i) Clearly, if $\alpha
\prec^\rightarrow \beta$, then $B_{\alpha,\beta} =
-(q - q^{-1}) q^{K_{\alpha,\beta}}$ since \\
$(e_{\beta'} \otimes e_{-\alpha'})(e_{\beta''} \otimes e_{-\alpha''})
\neq 0$ only if $\alpha' \prec^\leftarrow \beta'$, $\alpha''
\prec^\leftarrow \beta''$, and in this case $\alpha'+\alpha''
\prec^\leftarrow \beta' + \beta''$.  

(ii) We prove the lemma inductively.  If $|\alpha| = 1$, (ii) is
clear.  Otherwise, assume (ii) holds for $|\alpha| \leq p$.  We will prove
the result for $|\alpha| = p+1$.

Suppose $\alpha = e_i - e_{i+p+1}$, $\beta = e_j - e_{j+p+1}$, and
$\tau(\alpha_{i+k}) = \alpha_{j+p-k}$, $0 \leq k \leq p$.  Then, by
\eqref{jiprod}, we may write
\begin{multline*}
B_{\alpha, \beta} = -\text{sign}(\alpha, \beta) (q - q^{-1})
q^{K_{\alpha,\beta}} - \sum_{l = 1}^p A_{e_i - e_{i+l},
e_{j+p+1-l}-e_{j+p+1}} B_{e_{i+l}-e_{i+p+1},e_j-e_{j+p+1-l}} \\ =
-\sab \qh q^{\kab} - \sum_{l=1}^p \sab \qh^2
q^{p+1-2l+L_{\alpha,\beta}} \\ = -\sab \qh q^{L_{\alpha,\beta}}
\biggl[ q^{-p} + \sum_{l=1}^p \qh q^{p-2l+1} \biggr] = -\sab \qh
q^{p+L_{\alpha,\beta}}. \quad\square
\end{multline*}

Using $e_{-\alpha} \wedge_c e_\beta$ as in Example \ref{gcge}, \eqref{jdds} becomes
\begin{equation} \label{rjdf}
\bar R_J = R_s + (q - q^{-1}) \biggl[ \sum_{\alpha \prec^\rightarrow \beta}
e_{-\alpha} \wedge_{-\frac{1}{2}(\alpha, \beta)} e_{\beta} + \sum_{\alpha \prec^\leftarrow \beta} (-1)^{|\alpha|-1} e_{-\alpha} 
\wedge_{-\frac{1}{2}(\alpha, \beta)+|\alpha|-1} e_{\beta} \biggr].
\end{equation}

All that remains is to show equivalence of \eqref{rjdf} $\bar
R_{\text{GGS}}$.  First we write $\bar R_{\text{GGS}}$:

\begin{equation} \label{rggs}
\bar R_{\text{GGS}} = R_s + (q - q^{-1}) \sum_{\alpha \prec \beta} 
\text{sign}(\alpha, \beta) e_{-\alpha} \wedge_{-\text{sign}(\alpha, \beta)\epsilon_{\alpha,\beta}} e_{\beta}
\end{equation}

Combining \eqref{rggs} with Table 4.2, we obtain \eqref{rjdf}.
This proves that $\bar R_{\text{GGS}} = \bar R_J$ and hence
that $R_{\text{GGS}} = R_J$ in the case $\Gamma_1 \cap \Gamma_2
= \emptyset$.  The proof is finished.$\quad \square$

%In Appendix \ref{da} we prove the twist conjecture in the disjoint case
%by showing $R_J$ satisfies the QYBE (proceeding as in \cite{H}, \cite{ER}).

\subsection{The generalized disjoint case} \label{gds}

In fact, the results we have obtained extend easily to the {\it generalized
disjoint} case:

\begin{defe}  A triple $(\Gamma_1, \Gamma_2, \tau)$ is said to be
{\it generalized disjoint} if $\Gamma_1 = \bigcup_{i=1}^{m}
\Gamma_1^i$ where $\Gamma_1^i \perp \Gamma_1^j, i \neq j$ and
$\tau(\Gamma_1^i) \cap \Gamma_1 \subset \Gamma_1^{i+1}, i < m,$ while
$\tau \Gamma_1^m \subset \Gamma_2$.
\end{defe}

\begin{thm}\label{gdt}  For any generalized disjoint triple, $R_J = R_{\text{GGS}}$.
\end{thm}

{\it Proof.} Note first by \eqref{fk}, \eqref{fe} that $K_{\alpha,\beta}$ and
$\epsilon_{\alpha,\beta}$ are as given in Tables 4.1, 4.2,
respectively.  As in the disjoint case, we have the following main
property:

\begin{lemma} If $(e_{-\alpha} \o e_{\tau^x \alpha}) 
(e_{-\beta} \o e_{\tau^y \beta})_{ab} \neq 0$, for $\{a,b\} =
\{1,2\}$, then $(a,b)=(1,2)$, $x = y$, and $\tau^x$ reverses
orientation on $\alpha+\beta$.
\end{lemma}

{\it Proof.}  Suppose $(a,b)=(2,1)$.  Then $\{\tau^x \alpha, \beta\} \subset
\tilde \Gamma_1^i$ for some $i$.  But then, $\alpha \in \Gamma_1^j, \tau^y
\beta \in \Gamma_1^k$ for $j < i < k$, so $e_{-\alpha} e_{\tau^y \beta} = 0$,
a contradiction.  So, $(a,b)=(1,2)$.  Now, this implies that $\{\alpha, \beta\}
\subset \Gamma_1^i, \{\tau^x \alpha, \tau^y \alpha\} \subset \Gamma_1^j$
for some $i < j$.  Then $x=y=j-i.\quad \square$

Because of this fact, we may set $J^{ij}$ to be the $J$-matrix
corresponding to the disjoint triple $(\Gamma_1^i,\tau \Gamma_1^{j-1},
\tau^{j-i})$, and similarly define $\bar R_J^{ij},
R_{\text{GGS}}^{ij}$, and $A^{ij} = J^{ij} - 1$, so that $A^{ij} X
A^{kl} = 0$ whenever $(i,j) \neq (k,l)$. Thus, $J = \prod_{i < j}
J^{ij} = 1 + \sum_{i < j} A^{ij}$, and $A^{ij} X A^{kl} = 0$ whenever
$i \neq j$.  Hence, $\bar R_J = R_s + \sum_{i,j} \bar (R_J^{ij}-R_s) = 
R_s + \sum_{i,j} (\bar R_{\text{GGS}}^{ij}-R_s) = \bar R_{\text{GGS}}.
\quad\square$

In the appendix we prove the twist conjecture and hence the GGS
conjecture in the orthogonal generalized disjoint case (defined in
Section \ref{it}) by demonstrating that $R_J$ satisfies the QYBE.

%\begin{rem} \label{wto}  Using this proof for the generalized disjoint 
%case, one
%may define a {\it weakly} $\tau-orthogonal$ union to be an orthogonal
%union $(\Gamma_1,\Gamma_2,\tau) = \bigcup_i
%(\Gamma_1^{(i)},\Gamma_2^{(i)})$ where $\Gamma_2^{(i)} \cap
%\Gamma_1^{(j)} = \emptyset$ whenever $i < j$.  In this case, one
%may show that if $R_J = R_{\text{GGS}}$ for each subtriple in the union,
%the same holds for the union.
%\end{rem}

%In this case, one may show that any $\tau
%A triple that
%is a $\tau$-disjoint orthogonal union of two nonempty triples is
%called {\it decomposable}; otherwise it is {\it indecomposable}.  

\section{The Cremmer-Gervais triple}
\label{crgs}

In this section we prove $R_J = R_{\text{GGS}}$ in the Cremmer-Gervais
case, and hence the twist conjecture since the GGS conjecture is proved
in this case (for example, see \cite{H2}.)

\begin{thm}  For the Cremmer-Gervais triple 
$(\{\alpha_1,\ldots,\alpha_{n-1}\},\{\alpha_2,\ldots,
\alpha_n\}, \tau)$, \\
$\tau \alpha_i = \alpha_{i+1}$, one has $R_J = R_{\text{GGS}}$.
\end{thm}

{\it Proof.}  Note first from \eqref{fk}, \eqref{fe} that
$K_{\alpha,\beta} = \frac{1}{2}[\alpha \lessdot \beta] + [\alpha \ll
\beta]$, and $\epsilon_{\alpha,\beta} = -K_{\alpha,\beta}$.  Next,
note that $(e_{-\alpha} \o e_{\tau^x \alpha})(e_{-\beta} \o e_{\tau^y
\beta})_{ab} \neq 0$ for $\{a,b\} = \{1,2\}$ iff $(a,b) = (1,2)$ and
$y = x+|\alpha|+|\beta|$.  Hence, considering the product form
\eqref{jiprod}, we find that $J^{-1} = 1 + \sum_{\alpha \prec \beta}
q^{K_{\alpha,\beta}} (q^{-1} - q) \eab$.  Furthermore, if we set $B =
J^{-1} - 1, A = J - 1$, we see that 

\begin{multline} \label{cgrj}
\bar R_J =
J^{-1} R_s J_{21} = R_s + B + A_{21} + (q-1) A_{21}^{\lessdot}
+ (q-q^{-1}) P A_{21}^{\ll} \\ + (q-q^{-1}) P B_{21} +
(q-q^{-1}) P B_{21} A_{21} + B A_{21}.
\end{multline}

Since, in addition, $(e_{\beta} \o e_{-\alpha})_{ab} \in G$ for
$\alpha,\beta > 0, \{a,b\} = \{1,2\}$ iff $\alpha \prec \beta$, we may
infer that $(\bar R_J)_+ = (\bar R_{\text{GGS}})_+$.  Thus, it
suffices to show $(\bar R_J)_{-\beta,-\alpha} = (\bar
R_{\text{GGS}})_{-\beta,-\alpha}$ whenever $\alpha \prec \beta$
(equivalently, $\alpha < \beta$.)  We proceed inductively on the
index under $\tau$.  Suppose this is true for all $\alpha' \prec \beta'$
where $i_\tau(\beta'-\alpha') < i_\tau(\beta-\alpha)$.  
Since, with respect to the ordering in \eqref{jiprod}, $(e_{-\alpha} \o
e_{\beta})(e_{-\alpha'} \o e_{\beta'}) \neq 0$ only when $(\alpha,\beta)
< (\alpha',\beta')$, we may rewrite \eqref{cgrj} using our inductive
hypothesis as
\begin{multline} \label{icg}
(\bar R_J)_{-\beta,-\alpha} = (q-q^{-1}) \bigl[
q^{K_{\alpha,\beta}+[\alpha
\lessdot \beta]} + [\alpha \overline{<} \beta] q +
(P B_{21})_{-\beta,-\alpha} \\ + \bigl(\bar R_{\text{GGS}}
\sum_{\alpha \prec \beta} q^{K_{\alpha,\beta}} (e_{-\alpha} \o e_\beta)
\bigl)_{-\beta,-\alpha} \bigr].
\end{multline}

When we write out the sum above, we will add it in pairs corresponding
to the bijections $f,g$ defined in Section 3; this will make the sum
cancel.
%
%we may write $(\bar
%R_J)_{e_{i+k-j}-e_k,e_i - e_j} = \sum_{l = j+1}^i (q-q^{-1})
%q^{K_{e_j-e_l,e_{i+k-l}-e_{i+k-j}}} (R_J)_{e_{i+k-l}-e_k,e_i-e_l}$.
%We use this formula to inductively prove $(\bar R_J)_{-\beta,-\alpha}
%= (\bar R_{\text{GGS}})_{-\beta,-\alpha}$
%
We separately consider the cases $\alpha \overline{<} \beta,
\alpha \lessdot \beta, \alpha \ll \beta$.  

First suppose $\alpha \overline{<} \beta$.  Set $\alpha = e_j - e_i,
\beta = e_k - e_{i+k-j}$ where $j < k < i < i+k-j$.  Then
$(R_{\text{GGS}})_{-\beta, -\alpha} = q-q^{-1}$.  By \eqref{icg}, we
may write

\begin{multline}
(\bar R_J)_{-\beta,-\alpha} = (q-q^{-1}) \bigl[ 1 + q - q + \sum_{l = j+1}^k
\bigl((q-q^{-1})(e_{il} \o e_{k,i+k-l})q (e_{lj} \o e_{i+k-l,i+k-j}) \\ - q(q-q^{-1})(e_{i,i+l-j} \o e_{k,k-l+j})q^0(e_{i+l-j,j} \o e_{k-l+j,i+k-j}) \bigr)_{-\beta,-\alpha} \bigl] = q-q^{-1}.
\end{multline}

Now, suppose $\alpha \lessdot \beta$.  In this case, set $\alpha = e_j - e_i,
\beta = e_i - e_{2i-j}$.  Now, \eqref{icg} becomes

\begin{multline}
(\bar R_J)_{-\beta,-\alpha} = (q-q^{-1})\bigl[ q^{3/2} - q^{1/2}(q - q^{-1})
+ \sum_{l = i+1}^j \bigl( q^{-1/2}(q - q^{-1}) (e_{il} \o e_{i,2i-l}) q (e_{lj}
\o e_{2i-l,2i-j}) \\ - q^{1/2}(q-q^{-1}) (e_{i,2i-l} \o e_{i,l}) q^0 (e_{2i-l,j}
\o e_{l,2i-j}) \bigr)_{-\beta,-\alpha} \bigl] = q^{-1/2}(q - q^{-1}).
\end{multline}

Finally, suppose $\alpha \ll \beta$.  In this case, $\alpha = e_j - e_i,
\beta = e_k - e_{i+k-j}$ where $j < i < k < i+k-j$, and we write

\begin{multline}
(\bar R_J)_{-\beta,-\alpha} = (q-q^{-1}) \bigl[ q - (q - q^{-1}) + 
\sum_{l = j+1}^k
\bigl(q^{-1}(q-q^{-1})(e_{il} \o e_{k,i+k-l})q (e_{lj} \o e_{i+k-l,i+k-j}) \\ - (q-q^{-1})(e_{i,i+l-j} \o e_{k,k-l+j})q^0(e_{i+l-j,j} \o e_{k-l+j,i+k-j}) \bigr)_{-\beta,-\alpha} \bigl] = q^{-1}(q-q^{-1}).
\end{multline}

The proof is finished.\qs

%\begin{rem}  It is not difficult to extend this proof to the case
%when $\tau \alpha_i = \alpha_{i+k}, \Gamma_1 =
%\{\alpha_1,\ldots,\alpha_{n-k-1}\}, \Gamma_2 =
%\{\alpha_{k+1},\ldots,\alpha_{n-1}\},$ for any $k \in \N$.  A similar
%proof works for $\tau \alpha_i = \alpha_{i-k}$.  Combining this with
%Theorem \ref{gds} and some union arguments, one finds that $R_J =
%R_{\text{GGS}}$ whenever $\Gamma_1 = \bigcup \Gamma_1^{(i)}$, where
%each $\Gamma_1^{(i)}$ is connected and mutually orthogonal, and $\tau
%\Gamma_1^{(i)} \cap \Gamma_1^{(j)} = \emptyset$ for $i > j$.  In
%particular, this includes the case when $\tau \alpha_i = \alpha_j$
%implies $j > i$ for all $\alpha_i \in \Gamma_1$ (or $j < i$ for all
%$\alpha_i$).
%\end{rem}

%Now, by the
%product formula \eqref{jiprod}, we may consider $(R_J)_{-\beta,-\alpha}$
%as a sum $(q - q^{-1}) - q (q-q^{-1})^2$
% = q (q-q^{-1}) \sum_{l = j+1}^i 

%(i) $(\bar R_J)_{\alpha,\beta} = (\bar R_{\text{GGS}})_{\alpha,\beta} = q^{-1} - q$.  (ii) $(\bar R_J)_{-\beta,-\alpha} = (\bar R_{\text{GGS}})_{-\beta,-\alpha} = q - q^{-1}$.
%\end{lemma}

%{\it Proof.}  Inductively suppose that the result holds for all $\alpha

\section{Acknowledgements}

I would like to thank Pavel Etingof for introducing me to this problem
and advising me.  I would also like to thank the Harvard College
Research Program for their support.  Finally, I am indebted to
Gerstenhaber, Giaquinto, and Hodges for valuable discussions and for
sharing some unpublished results.

\appendix
\section{Proof that $R_J$ satisfies the QYBE in special cases}
\subsection{The disjoint case, by P.Etingof and T.Schedler}\label{da}
\vskip 12 pt

In this subsection we will prove Theorem \ref{tdg}.i by showing that $R_J$
satisfies the QYBE for disjoint triples (i.e. $\Gamma_1\cap
\Gamma_2=\emptyset$). We note that in the case when $\Gamma_1$ is orthogonal to
$\Gamma_2$, this was done (using the same method) by T.Hodges.  This 
is sufficient given Theorem \ref{dt}.

Let $U_\hbar(\mathfrak{sl}(n))$ be the quantum universal
enveloping algebra generated by $E_{\alpha_i}, F_{\alpha_i}$, $H_{\alpha_i}$,
$\alpha_i \in \Gamma,$ under the relations
\begin{gather} \label{ur1}
[H_{\alpha_i},E_{\alpha_j}] = (\alpha_i,\alpha_j) E_{\alpha_j},
[H_{\alpha_i},F_{\alpha_j}] = -(\alpha_i,\alpha_j) F_i,
[H_{\alpha_i},H_{\alpha_j}] = 0, \\ [E_{\alpha_i},F_{\alpha_j}] =
\delta_{ij} \frac{q^{H_{\alpha_i}} - q^{-H_{\alpha_i}}}{q - q^{-1}},\\
E_{\alpha_i}^2 E_{\alpha_{i \pm 1}} - (q + q^{-1})E_{\alpha_i}
E_{\alpha_{i \pm 1}} E_{\alpha_i} + E_{\alpha_{i \pm 1}}
E_{\alpha_i}^2 = 0,\\ F_{\alpha_i}^2 F_{\alpha_{i \pm 1}} - (q +
q^{-1})F_{\alpha_i} F_{\alpha_{i \pm 1}} F_{\alpha_i} + F_{\alpha_{i
\pm 1}} F_{\alpha_i}^2 = 0,
\end{gather}
where the coproduct, counit, and antipode are given by
\begin{gather} \label{ur2}
\Delta(E_{\alpha_i}) = E_{\alpha_i} \o q^{H_{\alpha_i}} + 1 \o E_{\alpha_i}, \Delta(F_{\alpha_i}) = F_{\alpha_i} \o 1 +
q^{H_{\alpha_i}} \o F_{\alpha_i},\\
\Delta(H_{\alpha_i}) = H_{\alpha_i} \o 1 + 1 \o H_{\alpha_i}, \epsilon(F_{\alpha_i})=\epsilon(E_{\alpha_i})=\epsilon(H_{\alpha_i}) = 0, \\
S(E_{\alpha_i}) = -E_{\alpha_i} q^{H_{\alpha_i}}, S(F_{\alpha_i}) = -q^{H_{\alpha_i}} F_{\alpha_i}, S(H_{\alpha_i}) = - H_{\alpha_i}.
\end{gather}

We will use the representation $\phi: U_\hbar(\mathfrak{sl}(n))
\rightarrow Mat_n(\C)$ by $\phi(E_{\alpha_i}) = e_{i,i+1},
\phi(F_{\alpha_i}) = e_{i+1,i}$, and $\phi(H_{\alpha_i}) = e_{ii} -
e_{i+1,i+1}$.

%Set $E_{\alpha_i} = E_i, F_{\alpha_i} = F_i, H_{\alpha_i} = H_i -
%H_{i+1}$, $1 \leq i \leq n$.  

Now, we recall the results of Hodges \cite{H} using the notation of
\cite{ER}. Fix a disjoint Belavin-Drinfeld triple
$(\Gamma_1,\Gamma_2,\tau)$.  Let $\h_i$ be the subpaces of $\h$
spanned by $e_{\alpha_k}, \alpha_k \in \Gamma_i$.  Let $U_i$ be the
Hopf subalgebras of $U_\hbar(\mathfrak{sl}(n))$ generated by
$E_{\alpha_k}, F_{\alpha_k}, H_{\alpha_k}, \alpha_k \in \Gamma_i$.
Define $f_{\tau}: U_1 \rightarrow U_2$ by $f_{\tau}(E_{\alpha_i}) =
E_{\tau \alpha_i}, f_{\tau} F_{\alpha_i} = F_{\tau \alpha_i}, f_\tau
H_{\alpha_i} = H_{\tau \alpha_i}$.  It is clear that $f_\tau$ is an
isomorphism of Hopf subalgebras. Let $g_\tau: \phi(U_1) \rightarrow
\phi(U_2)$ be the homomorphism descending from $f_\tau$.

Define $Z=(g_\tau \o 1) \Omega_{\h_1}$ where $\Omega_{\mathfrak{t}}$ denotes
the Casimir element of the usual bilinear form on a nondegenerate subspace 
$\mathfrak{t} \subset \h$. Let $T\in \h\o\h$ be a solution of the
following equations:
\begin{equation} \label{te}
(x\o 1,T)=(1\o\tau(x),T)=0,\quad (\tau(x)\o 1+1\o x,Z-T)=0.
\end{equation}
for any $x\in \h_1$. Define 
\begin{equation}
{\cal J}=q^T{\cal J}_0,\quad 
{\cal J}_0=1+\sum_{\alpha \prec \beta}(q-q^{-1})(-q^{-1})^{C_{\alpha}
(|\alpha|-1)}
e_{\beta}\o e_{-\alpha},
\end{equation}
where $C_\alpha=1$ if $\alpha \prec^{\leftarrow} \beta$ and $0$ if $\alpha \prec^{\rightarrow} \beta$. 

The following proposition can be deduced from the results of \cite{H}.

\begin{prop} \label{dtp} The element 
$R_{{\cal J}} \equiv {\cal J}^{-1}R_s{\cal J}_{21}\in Mat_n(\C)\o Mat_n(\C)$ 
satisfies the Hecke relation and the quantum Yang-Baxter equation. 
\end{prop} 
 
{\it Proof.} Define
\begin{equation} \label{jp}
{\cal J}'=q^{T-Z}(\tau\o 1)(\phi \o \phi)(\mathbb R)
\end{equation}
where $\mathbb R$ is the universal $R$-matrix of the Hopf subalgebra
$U_1$.  By Proposition 4.1 of \cite{ER}, $({\cal J}')^{-1}R_s({\cal
J}')_{21}$ satisfies the Hecke relation and the quantum Yang-Baxter
equation, so it is enough to show that ${\cal J}$ coincides with
${\cal J}'$.  This can be deduced from an explicit formula for the
universal $R$-matrix $\mathbb R$.  We use the formula given in
\cite{KhT} and evaluate it in the representation as follows:

\begin{equation} \label{er}
(\phi \o \phi)\bigl((f_\tau\o 1)(\mathbb R)\bigr)=
\bigl( 1+\sum_{\alpha \prec \beta}(q-q^{-1})(-q^{-1})^{[\alpha \prec^\leftarrow \beta]
(|\alpha|-1)}
e_{\beta}\o e_{-\alpha} \bigr) q^Z.
\end{equation}

Here we take our normal ordering to be given by left to right on the
Dynkin diagram (for simple roots).  The additional powers of $-q^{-1}$
in the reversing case appear because in \cite{KhT}, $E_{\alpha+\beta}
= E_{\alpha} E_{\beta} - q^{(\alpha,\beta)} E_{\beta} E_{\alpha},$ and
so $\phi(E_{\tau \alpha}) = (-q^{-1})^{|\alpha|-1} e_{\tau \alpha}$
when $\tau$ reverses order on $\alpha$, while $\phi(E_{\tau \alpha})
= e_{\tau \alpha}$ when $\tau$ preserves order on $\alpha$.

All that remains is to show $[q^Z, (\phi \o \phi)({\cal J}_0)] = 0.$  Observe
that $\sum_i e_{ii} \o e_{ii} = \Omega_h = \Omega_{h_1} + \Omega_{h_1^\perp}$
so that $[\Omega_{h_1}, (\phi \o \phi) (\R)] = 0$. This finishes the proof. \qs

We denote the usual inner product on $(\h_1 \oplus \h_2)$ by $I(x,y)$,
so that we may define the bilinear form $B(x,y) \equiv (x, \tau y)$
defined on $(\h_1 \oplus \h_2) \o \h_1$.  Similarly define $B^T(x,y) =
(\tau x, y)$ on $\h_1 \o (\h_1 \oplus \h_2)$.  Since $\h_1+\h_2$ has a
nondegenerate inner product, any element of $(\h_1\oplus \h_2)^{\o 2}$
can be regarded as a bilinear form on $\h_1\oplus \h_2$, by $X\mapsto
F_X(a,b)=(a\o b,X)$.  Thus such an element can be written as a 2 by 2
matrix whose ij-th entry is a form on $\h_i\o \h_j$. We will use such
notation below.
  
\begin{lemma} There exists a unique solution $T$ of the equations \eqref{te} 
in $(\h_1\oplus \h_2)^{\o 2}$. This solution has the form 
\begin{equation}
T=\left(\matrix 0&0\\ I+B&0\endmatrix\right).
\end{equation}
\end{lemma}

{\it Proof.} The proof is by a direct computation.$\quad \square$

Now let us compute $R_{{\cal J}}$. We have 
\begin{equation} \label{a1}
R_{{\cal J}}=
{\cal J}_0^{-1}q^{-T}R_sq^{T_{21}}({\cal J}_0)_{21}.
\end{equation}
Using the fact that 
\begin{equation} \label{s2cr}
q^XR_sq^{-X}=R_s, X\in S^2\h
\end{equation} 
one transforms \eqref{a1} to 
\begin{equation}
R_{{\cal J}}=
{\cal J}_0^{-1}q^{T_{21}}R_sq^{-T}({\cal J}_0)_{21}.
\end{equation}
It is clear that $T_{21}$ commutes with ${\cal J}_0$, since both components
commute separately,
thus
\begin{equation}
R_{{\cal J}}=
q^{T_{21}}{\cal J}_0^{-1}R_s({\cal J}_0)_{21}q^{-T}.
\end{equation}
Let $J$ be as defined in \eqref{taid} and \eqref{jd}, using Table
1.  Then $J=q^Y{\cal J}_0q^{-Y}$, where $Y=\frac{1}{2}\sum
e_{ii}\o e_{ii}$.  Thus, using \eqref{s2cr} again, we get
\begin{equation}
R_{{\cal J}}=
q^{T_{21}-Y} J^{-1}R_s J_{21}q^{-T+Y}.
\end{equation}
Setting $R_J = q^{r^0} J^{-1} R_s J_{21} q^{r^0}$, we get
\begin{equation}
R_{{\cal J}}=
q^{T_{21}-Y-r^0} R_J q^{-T+Y-r^0},
\end{equation}
%where $R$ is the R-matrix of the GGS conjecture, and $r^0$ is 
%as in Section 1.2. 
%This shows that 
%\begin{equation}
%R=q^{-T_{21}+Y+r^0}R_{{\cal J}}q^{T-Y+r^0}.
%\end{equation}
Write $Y$ as $Y_0+Y'$, where $Y_0\in (\h_1\oplus \h_2)^{\o 2}$, $Y'\in 
((\h_1\oplus \h_2)^\perp)^{\o 2}$. It is clear that $q^{Y'}$ commutes 
with $R_J$, so we obtain
\begin{equation}
R_{\cal J} = q^{T_{21}-Y_0-r^0}R_J q^{-T+Y_0-r^0}.
\end{equation}
It is easy to see that $Y_0=\frac{1}{2}
\left(\matrix I&I\\ I&I\endmatrix\right)$, 
and $r^0$ can be chosen to be the element of $(\h_1\oplus \h_2)^{\o 2}$ 
given by the matrix $\left(\matrix 0&\frac{1}{2}(I+B^T)\\ 
-\frac{1}{2}(I+B)& 0\endmatrix\right)$. Therefore, 
$R_{\cal J} = q^{-U} R_J q^U$,
where $U=\left(\matrix I/2&-B^T/2\\ -B/2&I/2\endmatrix\right)$.
 
Let $W=\left(\matrix I/2&B^T/2\\ B/2&I/2\endmatrix\right)$.  Then $W\in
S^2K$, where $K$ is the Lie algebra of symmetries of $R_J$, i.e. $K =
\{x \in Mat_n(\C) \mid [1 \o x + x \o 1,R_J] = 0\}$. Therefore,
$q^WR_Jq^{-W}=(D\o D)R_J(D^{-1}\o D^{-1})$ for a suitable diagonal matrix
$D$.  So in order to finish the proof of the twist conjecture,
it suffices to show that $q^{U+W}R_{{\cal J}}q^{-U-W}$ is a solution of
the Yang-Baxter equation.

Note that $U+W=\left(\matrix I&0\\ 0&I\endmatrix\right)$ and therefore
$U+W \in S^2(\h_2^\perp) \oplus S^2(\h_1^\perp)$.  Thus, 
the twist conjecture follows from the following lemma.

\begin{lemma} Let $X$ be an element of $S^2(\h_1\oplus\h_2)$ 
which is orthogonal to $\h_1\o \h_2$. Then $q^XR_Jq^{-X}
=(D\o D)J(D^{-1}\o D^{-1})$ for a suitable 
diagonal matrix $D$. 
\end{lemma}

{\it Proof.} We can write $X$ as $X_1+X_2$, where $X_1$ is orthogonal to 
$\h_1$ in both components and $X_2$ to $\h_2$. It suffices to prove the 
lemma for each of them separately. Let us do it for $X_1$, for $X_2$ the 
proof is analogous.   

Let $E_\alpha=e_{\tau \alpha}\o e_{-\alpha}$. 
It is clear that $q^{X_1}E_\alpha q^{-X_1}=\lambda(\alpha)E_\alpha$ 
for some eigenvalue $\lambda(\alpha)$. All we need to show is that 
$\lambda(\alpha+\beta)$, when it is defined, equals 
$\lambda(\alpha)\lambda(\beta)$. 

We may assume that $\alpha\ \lessdot \beta$, i.e. $\alpha$ is to the
left of $\beta$.  If $\tau$ reverses orientation on $\alpha+\beta$,
i.e. $\alpha+\beta \prec^{\leftarrow} \tau(\alpha+\beta)$, then the
statement is clear since $E_{\alpha+\beta}=E_\beta E_\alpha$. If not,
we have $E_{\alpha+\beta}^{t_2}=E_\alpha^{t_2}E_\beta^{t_2}$ ($t_2$ is
transposition in the second component).  Thus, using that the second
component of $X_1$ commutes with $e_\alpha, e_\beta$, we get
\begin{equation}
(q^{X_1}E_{\alpha+\beta}q^{-X_1})^{t_2}=
q^{X_1}E_\alpha^{t_2}E_\beta^{t_2}q^{-X_1}=
\lambda(\alpha)\lambda(\beta)E_\alpha^{t_2}E_\beta^{t_2}=
\lambda(\alpha)\lambda(\beta)E_{\alpha+\beta}^{t_2}.
\end{equation}
This implies the lemma.  $\quad \square$

\subsection{The orthogonal generalized disjoint case} \label{gda}
In this subsection we prove Theorem \ref{tdg}.ii by showing $R_J$ satisfies
the QYBE in the orthogonal generalized disjoint case, defined in Definition
\ref{gdd}.  This is sufficient due to Theorem \ref{gdt}.

  Set $\h_{kl} =
\text{Span}(H_{\alpha_i} \mid \alpha_i \in \Gamma_1^k, \tau^{l-1}
\alpha_i \in \Gamma_1)$.  We also consider $U_{kl}$, the Hopf
subalgebra of $U_\hbar(\mathfrak{sl}(n))$ generated by $E_{\alpha_i},
F_{\alpha_i}, H_{\alpha_i},$ for $\alpha_i \in \h_{kl}$, and
$U'_{kl}$, the Hopf subalgebra of $U_\hbar(\mathfrak{sl}(n))$ generated
by $E_{\tau^l \alpha_i}, F_{\tau^l \alpha_i}, H_{\tau^l \alpha_i}$.  
Define the
map $f_{\tau^l}: U_{kl} \rightarrow U'_{kl}$ by 
$f_{\tau^l}(E_{\alpha_i}) = E_{\tau^l \alpha_i},
f_{\tau^l} F_{\alpha_i} = F_{\tau^l \alpha_i}, f_{\tau^l} H_{\alpha_i}
= H_{\tau^l \alpha_i}$.  
%Set $\h'_{kl} = \text{Span}(H_{\tau^l
%\alpha_i} \mid H_{\alpha_i} \in \h_{kl})$. 
Then it is easy to check
that $f_{\tau^l}$ is a Hopf algebra isomorphism.
% $U_{kl}
%\rightarrow U'_{kl}$ where $U'_{kl}$ is generated by
%$E_{\alpha_i},F_{\alpha_i},H_{\alpha_i}, \alpha_i \in \h'_{kl}$. Also,
%define the map $g_\tau: 
Also, define the map $g_{\tau^l}: \phi(U_{kl}) \rightarrow
\phi(U'_{kl})$ to be the homorphism descending from $f_{\tau^l}$.

Now, define $\R^{kl} = (f_{\tau^l} \o 1)(\R')^{kl}$ where $(\R')^{kl}$ is
the universal $R$-matrix of $U_{kl}$.  Now, we define the twist ${\cal
J}$ by ${\cal J} = \prod_{i = 1}^{m} \prod_{j = 1}^{m+1-i} \R^{j,j+i}$.
Recall from \cite{ER} the definition of a twist:

\begin{defe} ${\cal J} \in U_q(\mathfrak{sl}(n))$ is a twist if ${\cal J} \equiv 1 \pmod \hbar, (\epsilon \o 1) {\cal J} = (1 \o \epsilon) {\cal J} = 1,$ and $(\Delta \o 1) ({\cal J}) {\cal J}_{12} = (1 \o \Delta) ({\cal J}) {\cal J}_{23}$.  In addition, a twist is triangular
if ${\cal J} \in U_{\geq 0} \o U_{\leq 0}$, where $U_{\geq 0}, U_{\leq 0} 
\subset U_\hbar(\mathfrak{sl}(n))$ are the Hopf subalgebras generated by the
$E_{\alpha_i}, H_{\alpha_i}$, and by the $F_{\alpha_i}, H_{\alpha_i}$,
respectively.
\end{defe}

\begin{thm}  ${\cal J}$ is a triangular twist.
\end{thm}

{\it Proof.}  It is obvious that ${\cal J} \equiv 1 \pmod \hbar$,
$(\epsilon \o 1) {\cal J} = (1 \o \epsilon) {\cal J} = 1$, and that $\cal
J$ is triangular. It thus suffices to prove $(\Delta \o 1) ({\cal J})
{\cal J}_{23} = (1 \o \Delta) ({\cal J}) {\cal J}_{12}$.  By
construction, one easily sees that $[\R^{kl}_{ab}, \R^{k' l'}_{a' b'}] =
0$ for $1 \leq a < b \leq 3, 1 \leq k < l \leq m+1$ and similarly for
$a',b',k',l'$ in the event that $\{(a,k),(b,l)\} \cap
\{(a',k'),(b',l')\} = \emptyset$.  Furthermore, since $\R$ satisfies
the QYBE, where $\R$ is the universal $R$-matrix of any Hopf subalgebra
of $U_h(\mathfrak{sl}(n))$, we find $\R^{ij}_{12} \R^{ik}_{13}
\R^{jk}_{23} = \R^{jk}_{23} \R^{ik}_{13} \R^{ij}_{12}$ for $i < j < k$.
Thus, the theorem follows from the following combinatorial lemma:
 
\begin{lemma} Let G be a semigroup generated by the set $T=\{r_{ij}^{ab}, 1 \leq i < j \leq 3, 1 \leq a < b \leq n\}$ for
a given $n \in Z^+$.
Consider the relations
\begin{gather} \label{gr1}
r_{ij}^{ab} r_{kl}^{cd} = r_{kl}^{cd} r_{ij}^{ab}, \quad \{(i,a),(j,b)\} \cap
\{(k,c), (l,d)\} = \emptyset \\ \label{gr2}
r_{12}^{ab} r_{13}^{ac} r_{23}^{bc} = r_{23}^{bc} r_{13}^{ac} r_{12}^{ab},
\quad 1 \leq a < b < c \leq n.
\end{gather}
  Then, if $G$ satisfies relations \eqref{gr1}, \eqref{gr2}, it also
  satisfies the relation
\begin{equation} \label{gr3}
\biggl( \prod_{i=1}^{n-1} \prod_{j=1}^{n-i} r_{13}^{j,i+j} r_{12}^{j,i+j} \biggr)
\prod_{i=1}^{n-1} \prod_{j=1}^{n-i} r_{23}^{j,i+j} = 
\biggl( \prod_{i=1}^{n-1} \prod_{j=1}^{n-i} r_{13}^{j,i+j} r_{23}^{j,i+j} \biggr)
\prod_{i=1}^{n-1} \prod_{j=1}^{n-i} r_{12}^{j,i+j}.
\end{equation}
\end{lemma}
{\it Proof.}  Let $F$ be the free semigroup generated by the set $T$
as above, and let $Y$ be the set of relations \eqref{gr1},\eqref{gr2}.
Note that in each side of \eqref{gr3}, every generator of $G$ appears
exactly once.  Let $H \subset F$ be the set of such elements of $F$,
so that for any element $X \in H$, we can say $r_1 <_X r_2$ for
generators $r_1, r_2$ if $r_1$ appears to the left of $r_2$ in $X$.
Let $L,R$ denote the left and right hand sides of \eqref{gr3},
respectively, considered as elements of $H$.  Note that both $L$ and
$R$ satisfy
\begin{equation} \label{pr1}
r_{13}^{ab} <_X r_{12}^{ac}, r_{13}^{ef} <_X r_{23}^{df}, \quad b \leq c, d \leq e
\end{equation} 
when one replaces $X$ with $L,R$. Again replacing $X$ with $L,R$ we find that
\begin{equation} \label{pr2}
r_{12}^{ab} <_X r_{13}^{ac} <_X r_{23}^{bc}\quad \text{or}\quad r_{23}^{bc} <_X r_{13}^{ac} <_X r_{12}^{ab}, \quad 1 \leq a < b < c \leq n.
\end{equation}
Finally, consider the property
\begin{equation} \label{pr3}
r_{ij}^{ab} <_X r_{ij}^{cd} \Leftrightarrow r_{ij}^{ab} <_L r_{ij}^{cd}, \text{whenever}\ a = c\ \text{or}\ b=d.
\end{equation}
Obviously $L$ satisfies \eqref{pr3}, and it is easy to see that $R$
does as well.  Let $K \subset H$ be the set of elements of $H$
satisfying \eqref{pr1}, \eqref{pr2}, \eqref{pr3}.  Define a function
$f: K \rightarrow \Z^+$ by $f(X) = |\{(a,b,c) \in \Z^3| 1 \leq a < b <
c \leq n, r_{12}^{ab} <_X r_{13}^{ac} <_X r_{23}^{bc}\}|.$ We claim
that $|K/Y| = 1$, that is, the image of $K$ is just one element in the
natural map $F \rightarrow F/Y=G$.  To show this, first note that if
$X \in K$ satisfies $f(X) = 0$, then for all pairs of generators $r_1,
r_2 \in T$ that do not appear in \eqref{gr1}, $r_1 <_X r_2 \Leftrightarrow
r_1 <_L r_2$.  This follows from \eqref{pr1}, \eqref{pr2}, and
$f(X)=0$.  Thus, $f(X)=0 \Rightarrow X \equiv L \pmod Y$.  Now, we
show that given $f(X) = m$, $\exists Z \in H$ such that $X \equiv Z
\pmod Y$ and $f(Z) = m-1$, whenever $m \in \Z^+$.  Indeed, take a
triple $(a,b,c)$ so that $r_{23}^{bc} <_X r_{13}^{ac} <_X
r_{12}^{ab}$, and so that for any other triple $(a',b',c')$ satisfying
this property, $b'-a' \geq b-a$, and if $a = a', b = b'$, then $c' <
c$.  Clearly there exists such a triple.  

Now, consider terms $r \in T$ with
$r_{23}^{bc} <_X r <_X r_{13}^{ac}$.  We will consider all cases when
$r$ does not commute with $r_{23}^{bc}$ or $r_{13}^{ac}$ under
\eqref{gr1}.  There are two possibilities: (i)
$r=r_{23}^{dc}$. Now, \eqref{pr3} implies $d < b$, and \eqref{pr1} implies
$a < d$.  But this contradicts our choice of $a,b,c$.  (ii) $r =
r_{13}^{dc}$. Here \eqref{pr3} implies $d>a$, and \eqref{pr1} implies
$d<b$.  This contradicts our choice of $a,b,c$. Hence, $X \equiv X'
\pmod Y$ where $r_{23}^{bc}$ appears next to $r_{13}^{ac}$, and
$f(X')=m$.

Now, consider terms $r \in T$ with $r_{13}^{ac} <_{X'} r <_{X'} r_{12}^{ab}$
that don't commute with $r_{12}^{ab}$.  Contradictions will be made
with our choice of $a,b,c$.  There are four possibilities. (i)
$r=r_{12}^{ad}$: \eqref{pr3} implies $d < b$, contradiction.  (ii)
$r=r_{12}^{db}$: \eqref{pr3} implies $a < d$, contradiction.
(iii) $r=r_{13}^{ad}$: \eqref{pr3} implies $c < d$, contradiction. (iv) $r=r_{23}^{bd}$:  \eqref{pr3} implies $c < d$, contradiction.  Hence, $X' \equiv X'' \pmod Y$ where the terms $r_{23}^{bc},r_{13}^{ac},r_{12}^{ab}$ appear consecutively, and $X'' \equiv Z \pmod Y$ where Z differs from $X''$ by swapping $r_{23}^{bc},r_{12}^{ab}$.  Since $f(Z)=m-1$, the proof is finished. $\quad\square$

\begin{corr} The element ${\cal J}_{21}^{-1} {\cal R} {\cal J}$ satisfies the QYBE,
where $\cal R$ is the universal $R$-matrix for $U_h(\mathfrak{sl}(n))$.  Hence, 
$(\phi \o \phi)({\cal J})^{-1} R_s (\phi \o \phi)({\cal J})_{21}$ satisfies the QYBE.
\end{corr}

{\it Proof.} This is clear since $q^{-\frac{1}{n}} (R_s)_{21} =
(\phi \o \phi)({\cal R})$ (which is an easy consequence of the formula in
\cite{KhT}).\qs

\begin{prop} For a suitable $r^0$ and some $x \in S^2 \h$, 
$J q^{-r^0} = q^x (\phi \o \phi)({\cal J})$.
\end{prop}

{\it Proof.}  First, we explicitly evaluate $(\phi \o \phi) (\R^{kl})$ using the
formulas in \cite{KhT}.  Indeed, we may write 
\begin{equation}
(\phi \o \phi)(\R^{kl}) = \sum_{\alpha \in \tilde
\Gamma_k \cap \tau^{-l} \tilde \Gamma} (-q^{-1})^{[\alpha
\prec^{\leftarrow} \beta](|\alpha| - 1)} e_{\tau^l \alpha} \o
e_{-\alpha} q^{Z^{kl}},
\end{equation}
where $Z^{kl} = (g_{\tau^l} \o 1) \Omega_{\h_{kl}}$ where $\Omega_{\h_{kl}}$ is the Casimir
element for the space $\h_{kl}$ with the usual bilinear form.
The coefficients $(-q^{-1})^{[\alpha
\prec^{\leftarrow} \beta](|\alpha| - 1)}$ follow easily from the fact
that $E_{\alpha+\beta} = E_{\alpha} E_{\beta} - q^{(\alpha,\beta)}
E_{\beta} E_{\alpha}$ (as defined in \cite{KhT}) upon evaluation in
$Mat_n(\C) \o Mat_n(\C)$, since different terms vanish in the
reversing and non-reversing cases. Note that $(\phi \o \phi) (\R^{kl})_{\alpha,
\tau^l \alpha} = J_{\alpha, \tau^\alpha}$ for $\alpha \in \tilde
\Gamma_k \cap \tau^{-l} \tilde \Gamma$.  It remains to reconcile 
the extra terms, $q^{-r^0}$ and the $q^{Z^{kl}}$.

As in the previous subsection, one sees that $\Omega_{\h} =
\Omega_{\h_{kl}^\perp} + \Omega_{\h_{kl}}$ where $\h_{kl}^{\perp}$ is
the orthogonal complement to $\h_{kl}$ in $\h$.  Since $\Omega_{\h} =
\sum_i e_{ii} \o e_{ii}$, we have $[\Omega_{\h}, (\phi \o \phi)
(\R')^{kl}] = 0$, and it follows that $[q^{Z^{kl}}, (\phi \o
\phi)(\R^{kl})] = 0.$ By orthogonality, $[q^{Z^{kl}},(\phi \o
\phi)(\R^{k'l'})] = 0$ for $(k',l') \neq (k,l)$.  Now, set $Y =
\sum_{k,l} Z^{kl}$. Then, ${\cal J} = J q^Y$.  Now, it is clear that
$(\alpha \o \beta, Y) = (\alpha, \tau^l \beta)$ for $\alpha \in
\Gamma_1^{j} \cup \tau \Gamma_1^{j-1}, \beta \in \Gamma_1^{j+l} \cup
\tau \Gamma_1^{j+l-1}$ where we again take $\Gamma_1^0 =
\Gamma_1^{m+1} = \emptyset$.  Then, we may take $r^0 =
\frac{1}{2}(Y_{21} - Y)$, so that $r^0 + Y \in S^2 \h$.  Furthermore,
as we have seen, $[q^{Y}, {\cal J}]=0$, and it is clear that
$[q^{Y_{21}}, {\cal J}] = 0$ by orthogonality. The proof is
finished.\qs

\begin{corr}  $R_J$ satisfies the QYBE. Theorem \ref{tdg} is proved.\end{corr}
{\it Proof.}  Clear.\qs

\section{Proof of Giaquinto's formula \eqref{gcgr} for 
$R_{\text{GGS}}$ in the generalized Cremmer-Gervais case}

In this section we explicitly compute $R_{\text{GGS}}$ for generalized
Cremmer-Gervais triples, the only triples satisfying $|\Gamma_1| + 1 =
|\Gamma|$ (omitting only one root).  These are precisely the cases
where $r^0$ is unique if its first component has trace zero ($r^0 \in
\wedge^2 \h'$ where $\h' \subset \h$ is the subspace of diagonal
matrices of trace zero.)  First we summarize the results given in
Example \ref{gcge} as proved in \cite{GG}.  Let $\Res$ give the
residue mod $n$ in $\{1,\ldots,n\}$.  Take the triple indexed by
$(n,m)$, where $\tau \alpha_i = \alpha_{\Res(i+m)}$.  The unique $r^0$
whose first component has trace zero is given by $(r^0)_{ii}^{ii} = 0,
(r^0)^{ij}_{ij} = \frac{1}{2} - \Res\bigl(\frac{j-i}{m}).$ Then, the
only difficulty is in computing $q^{r^0} \tilde a q^{r^0}$, so here we
use \eqref{fe} to prove $q^{r^0} \tilde a q^{r^0} = \sum_{\alpha \prec
\beta} e_{-\alpha} \wedge_{\frac{-2O(\alpha,\beta)}{n}} e_\beta$.

Clearly we have 
\begin{equation}
q^{r^0} \tilde a q^{r^0} = \sum_{\alpha \prec \beta}
e_{-\alpha} \wedge_{\epsilon_{\alpha,\beta} + r(\alpha,\beta)}
e_{\beta},
\end{equation}
where $r(e_j - e_i, e_k - e_{i+k-j}) =
(r^0)_{j,i+k-j}^{j,i+k-j} + (r^0)_{i,k}^{i,k}, j < i$, 
since $\text{sign}(\alpha,\beta) = 1$ for all $\alpha \prec \beta$.
Thus, it suffices to show $\text{sign}(\alpha,\beta)
\epsilon_{\alpha,\beta} + r(\alpha,\beta) =
-\frac{2O(\alpha,\beta)}{n}$ where, as before, $O(\alpha,\beta) = m$
when $\tau^m \alpha = \beta$.  

One sees that 
\begin{multline}
(r^0)_{j,i+k-j}^{j,i+k-j} + (r^0)_{i,k}^{i,k} = 1 - \frac{1}{2}[2j =
i+k] - \frac{1}{2}[i = k] -
\frac{1}{n}\text{Res}\bigl(\frac{i+k-2j}{m}\bigr) -
\frac{1}{n}\text{Res} \bigl(\frac{k-i}{m}\bigr) \\ = 1 -
\frac{1}{2}[2j = i+k] - \frac{1}{2}[i = k] - \frac{2}{n}
\text{Res}\bigl(\frac{k-j}{m}\bigr) + M_{i,j,k},
\end{multline}
where
\begin{equation}
M_{i,j,k} = [2j \neq i+k][\text{Res}\bigl(\frac{k-j}{m}\bigr) >
\text{Res}\bigl(\frac{i+k-2j}{m})] -
[i \neq k][\text{Res}\bigl(\frac{k-j}{m}\bigr) <
\text{Res}\bigl(\frac{k-i}{m}\bigr)].
\end{equation}
 Thus, since
$\text{Res}\bigl(\frac{k-j}{m}\bigr) = O(\alpha,\beta)$, it suffices
to show $1 + M_{i,j,k} + \epsilon_{\alpha,\beta} = 0$.  In the case
$\alpha \lessdot \beta$ or $\beta \lessdot \alpha$, it is easy to see
that $\epsilon_{\alpha,\beta} = M_{i,j,k} = -\frac{1}{2}$.  Otherwise,
we may write
\begin{equation}
1+M_{i,j,k} = [\text{Res}\bigl(\frac{k-j}{m}\bigr) >
\text{Res}\bigl(\frac{i+k-2j}{m})] +
[\text{Res}\bigl(\frac{k-j}{m}\bigr) >
\text{Res}\bigl(\frac{k-i}{m}\bigr)],
\end{equation}
and it is easy to see that $[\text{Res}\bigl(\frac{k-j}{m}\bigr) >
\text{Res}\bigl(\frac{i+k-2j}{m})] = [\exists \gamma, \alpha \prec
\gamma \prec \beta, \alpha \gtrdot \gamma]$ while \\
$[\text{Res}\bigl(\frac{k-j}{m}\bigr) >
\text{Res}\bigl(\frac{k-i}{m}\bigr)] = [\exists \gamma, \alpha \prec
\gamma \prec \beta, \alpha \lessdot \gamma]$.  This finishes the proof. \qs

%\end{singlespace}
\end{document}